\documentclass[12pt,a4paper]{amsart}
\usepackage{amssymb,amsxtra}
\usepackage[english,french]{babel}

\theoremstyle{definition}

\theoremstyle{remark}

\theoremstyle{plain}

\theoremstyle{remark}

\newtheorem*{example}{Example}
\numberwithin{equation}{section}

\begin{document}

\title{ Polynomial equations 
 and  rank of matrices over $\mathbb{F}_{2} $ related to  persymmetric matrices}
\author{Jorgen~Cherly}
\address{D\'epartement de Math\'ematiques, Universit\'e de
    Brest, 29238 Brest cedex~3, France}
\email{Jorgen.Cherly@univ-brest.fr}
\email{andersen69@wanadoo.fr}

\maketitle

\begin{abstract}
Dans cet article nous illustrons par quelques exemples la connexion entre le nombre de solutions 
d'\' equations  polyn\^omiales satisfaisant des conditions de degr\' es et le nombre de rang i matrices rattach\' ees  aux
matrices  persym\' etriques
 \end{abstract}

\selectlanguage{english}

\begin{abstract}
In this paper we illustrate  by some examples the connection between the number of solutions of polynomial
equations satisfying degree conditions and the number of rank i matrices related to persymmetric matrices.

 \end{abstract}

\newpage

\maketitle 
\newpage
\tableofcontents
\newpage
\section{Notations}
\label{sec 1}
A matrix  $[\alpha_{i,j}] $ is persymmetric if  $\alpha _{i,j} = \alpha_{r,s}$   for   i+j=r+s.\\

  \vspace{0.1 cm}
We denote by $  \Gamma_{i}^{s\times k}  $  the number of rank i persymmetric $s\times k $ matrices over  $\mathbb{F}_{2}$ 
of the form :
  $$ \left ( \begin{array} {cccccc}
\alpha _{1} & \alpha _{2} & \alpha _{3} &  \ldots & \alpha _{k-1}  &  \alpha _{k} \\
\alpha _{2 } & \alpha _{3} & \alpha _{4}&  \ldots  &  \alpha _{k} &  \alpha _{k+1} \\
\vdots & \vdots & \vdots   &  \vdots  & \vdots  & \vdots \\
\vdots & \vdots & \vdots    &  \vdots & \vdots & \vdots \\
\alpha _{s-1} & \alpha _{s} & \alpha _{s+1} & \ldots  &  \alpha _{k+s-3} &  \alpha _{k+s-2}  \\
\alpha _{s} & \alpha _{s+1} & \alpha _{s+2} & \ldots  &  \alpha _{k+s-2} &  \alpha _{k+s-1}  \\
\end{array}  \right). $$  \vspace{0.1 cm}\\
We denote by  $ \Gamma_{i}^{\left[s\atop s+m \right]\times k} $
 the number of rank i double persymmetric $(2s+m) \times k $ matrices over  $\mathbb{F}_{2}$ 
of the form:\vspace{0.1 cm}
 $$   \left ( \begin{array} {cccccc}
\alpha _{1} & \alpha _{2} & \alpha _{3} &  \ldots & \alpha _{k-1}  &  \alpha _{k} \\
\alpha _{2 } & \alpha _{3} & \alpha _{4}&  \ldots  &  \alpha _{k} &  \alpha _{k+1} \\
\vdots & \vdots & \vdots    &  \vdots & \vdots  &  \vdots \\
\alpha _{s-1} & \alpha _{s} & \alpha _{s +1} & \ldots  &  \alpha _{s+k-3} &  \alpha _{s+k-2}  \\
 \alpha _{s} & \alpha _{s+1} & \alpha _{s +2} & \ldots  &  \alpha _{s+k-2} &  \alpha _{s+k-1} \\
 \hline \\
\beta  _{1} & \beta  _{2} & \beta  _{3} & \ldots  &  \beta_{k-1} &  \beta _{k}  \\
\beta  _{2} & \beta  _{3} & \beta  _{4} & \ldots  &  \beta_{k} &  \beta _{k+1}  \\
\vdots & \vdots & \vdots    &  \vdots & \vdots  &  \vdots \\
\beta  _{m+1} & \beta  _{m+2} & \beta  _{m+3} & \ldots  &  \beta_{k+m-1} &  \beta _{k+m}  \\
\vdots & \vdots & \vdots    &  \vdots & \vdots  &  \vdots \\
\beta  _{s+m-1} & \beta  _{s+m} & \beta  _{s+m+1} & \ldots  &  \beta_{s+m+k-3} &  \beta _{s+m+k-2}  \\
  \beta  _{s+m} & \beta  _{s+m+1} & \beta  _{s+m+2} & \ldots  &  \beta_{s+m+k-2} &  \beta _{s+m+k-1}
\end{array}  \right). $$ \vspace{0.1 cm}\\
We denote by  $ \Gamma_{i}^{\left[s\atop{ s+m\atop s+m+l} \right]\times k} $
 the number of rank i double persymmetric $(3s+2m+l) \times k $ matrices over  $\mathbb{F}_{2}$ 
of the form:\vspace{0.1 cm}
 $$   \left ( \begin{array} {cccccc}
\alpha _{1} & \alpha _{2}  &  \ldots & \alpha _{k-1}  &  \alpha _{k} \\
\alpha _{2 } & \alpha _{3} &  \ldots  &  \alpha _{k} &  \alpha _{k+1} \\
\vdots & \vdots & \vdots    & \vdots  &  \vdots \\
\alpha _{s-1} & \alpha _{s} & \ldots  &  \alpha _{s+k-3} &  \alpha _{s+k-2}  \\
 \alpha  _{s } & \alpha  _{s +1} & \ldots & \alpha  _{s +k-2}& \alpha  _{s +k-1}\\
\hline \\
\beta  _{1} & \beta  _{2}  & \ldots  &  \beta_{k-1} &  \beta _{k}  \\
\beta  _{2} & \beta  _{3}  & \ldots  &  \beta_{k} &  \beta _{k+1}  \\
\vdots & \vdots    &  \vdots & \vdots  &  \vdots \\
\beta  _{m+1} & \beta  _{m+2}  & \ldots  &  \beta_{k+m-1} &  \beta _{k+m}  \\
\vdots & \vdots    &  \vdots & \vdots  &  \vdots \\
\beta  _{s+m-1} & \beta  _{s+m}  & \ldots  &  \beta_{s+m+k-3} &  \beta _{s+m+k-2}  \\
 \beta _{s+m} & \beta _{s+m+1} & \ldots & \beta _{s+m+k-2} & \beta _{s+m+k-1}\\
\hline \\
\gamma  _{1} & \gamma   _{2}  & \ldots  & \gamma  _{k-1} &  \gamma  _{k}  \\
\gamma  _{2} & \gamma  _{3}  & \ldots  & \gamma  _{k} &  \gamma  _{k+1}  \\
\vdots & \vdots    &  \vdots & \vdots  &  \vdots \\
 \gamma  _{m+1} &  \gamma _{m+2}  & \ldots  & \gamma _{k+m-1} &  \gamma  _{k+m}  \\
\vdots & \vdots   &  \vdots & \vdots  &  \vdots \\
 \gamma  _{s+m-1} & \gamma  _{s+m}  & \ldots  & \gamma  _{s+m+k-3} &  \gamma  _{s+m+k-2}  \\
  \gamma  _{s+m} & \gamma  _{s+m+1}  & \ldots  & \gamma  _{s+m+k-2} &  \gamma  _{s+m+k-1}\\
\gamma  _{s+m+1} & \gamma  _{s+m+2}  & \ldots  & \gamma  _{s+m+k-1} &  \gamma  _{s+m+k}\\
\vdots & \vdots   &  \vdots & \vdots  &  \vdots \\
\gamma  _{s+m+l} & \gamma  _{s+m+l+1}  & \ldots  & \gamma  _{s+m+l+k-2} &  \gamma  _{s+m+l+k-1}
\end{array}  \right). $$ 

\newpage
\section{Polynomial equations related to persymmetric matrices}
\label{sec 2 }

\begin{example}

The number   $  \Gamma_{i}^{5\times 5}  $ of persymmetric $5\times 5$ matrices  over  $\mathbb{F}_{2}$   of rank  i of the form \\
 $$   \left ( \begin{array} {ccccc}
\alpha  _{1} & \alpha  _{2}  &   \alpha_{3} &  \alpha _{4} &  \alpha _{5}  \\
\alpha  _{2} & \alpha  _{3}  &   \alpha_{4} &  \alpha _{5} &  \alpha _{6}  \\
\alpha  _{3} & \alpha  _{4}  &   \alpha_{5} &  \alpha _{6} &  \alpha _{7}  \\
\alpha  _{4} & \alpha  _{5}  &   \alpha_{6} &  \alpha _{7} &  \alpha _{8}  \\
\alpha  _{5} & \alpha  _{6}  &   \alpha_{7} &  \alpha _{8} &  \alpha _{9}  
 \end{array}  \right). $$ 
  is equal to \\
  \begin{equation*}
\left\{\begin{array}{ccc}
             1 & if & i= 0, \\
              3 & if & i= 1, \\
              12 & if & i= 2, \\
                48 & if & i= 3, \\
                 192 & if & i= 4, \\
                  256 & if & i= 5.
             \end{array}\right.\
\end{equation*}

The general formula is given by  D.E. Daykin [1] and in [2] in the case q = 2.\\

\textbf{Application:}
 The number $ R_{q} $ of solutions
 $(Y_1,Z_1, \ldots,Y_q,Z_q) \in \big(\mathbb{F}_{2}[T] \big)^{2q} $  of the polynomial equation
                        $$ Y_1Z_1 +  Y_2Z_2 + \ldots + Y_qZ_q = 0  $$ \\
  satisfying the degree conditions \\
                   $$  degY_j \leq 4 , \quad degZ_j \leq 4 \quad for  1\leq j \leq q, $$ \\
   are connected with the numbers  $ \Gamma_{i} $ of persymmetric 
    $5\times 5$ matrices over  $\mathbb{F}_{2} $
     of rank i, in the way that the number $ R_{q} $ can be written as a linear combination of the $\Gamma_{i}.$\\
  More precisely:\\
    \begin{align*}
& R_{q} 
  = 2^{10q-9}\sum_{i=0}^{5}\Gamma_{i}2^{-qi}  \\
 &  = \left\{\begin{array}{ccc}
   63 & if & q=1,\\
   8704 & if & q=2, \\
   2^{10q-9}\left[1+ 3\cdot2^{-q} +12\cdot2^{-2q} + 48\cdot2^{-3q} +192\cdot2^{-4q} +256\cdot2^{-5q} \right]
 & if & 3\leq q.
  \end{array}\right. 
  \end{align*}
\end{example}

\begin{example}

The number   $ \Gamma_{i} = \Gamma_{i}^{2\times 2}  $ of persymmetric $2\times 2$ matrices  over  $\mathbb{F}_{2}$   of rank  i of the form \\
 $$   \left ( \begin{array} {ccccc}
\alpha  _{1} & \alpha  _{2}    \\
\alpha  _{2} & \alpha  _{3}   
 \end{array}  \right). $$ 
  is equal to \\
  \begin{equation*}
\left\{\begin{array}{ccc}
             1 & if & i= 0, \\
              3 & if & i= 1, \\
              4& if & i= 2, \\
               \end{array}\right.\
\end{equation*}

\textbf{Application:}
 The number $ R_{q} $ of solutions
 $(Y_1,Z_1, \ldots,Y_q,Z_q) \in \big(\mathbb{F}_{2}[T] \big)^{2q} $  of the polynomial equation
                        $$ Y_1Z_1 +  Y_2Z_2 + \ldots + Y_qZ_q = 0  $$ \\
  satisfying the degree conditions \\
                   $$  degY_j \leq 1 , \quad degZ_j \leq 1 \quad for  1\leq j \leq q, $$ \\
   is connected with the numbers  $ \Gamma_{i} $ of persymmetric 
    $2\times 2$ matrices over  $\mathbb{F}_{2} $
     of rank i, in the way that the number $ R_{q} $ can be written as a linear combination of the $\Gamma_{i}.$\\
  More precisely:\\
    \begin{align*}
& R_{q} 
  = 2^{4q-3}\sum_{i=0}^{2}\Gamma_{i}2^{-qi}  \\
 &  = \left\{\begin{array}{ccc}
   7 & if & q=1,\\
   64 & if & q=2, \\
   2^{4q-3}\left[1+ 3\cdot2^{-q} +4\cdot2^{-2q}  \right]
 & if & 3\leq q.
  \end{array}\right. 
  \end{align*}
\end{example}

\begin{example}
The fraction of square persymmetric matrices which are invertible is equal to $ {1\over 2} $
\end{example}

\begin{example}

 The number  $ \Gamma _{i}^{\Big[\substack{5 \\ 1 }\Big] \times 5}$ of rank i matrices of the form \\
  $$  \left ( \begin{array} {ccccc}
\alpha  _{1} & \alpha  _{2}  &   \alpha_{3} &  \alpha _{4} &  \alpha _{5}  \\
\alpha  _{2} & \alpha  _{3}  &   \alpha_{4} &  \alpha _{5} &  \alpha _{6}  \\
\alpha  _{3} & \alpha  _{4}  &   \alpha_{5} &  \alpha _{6} &  \alpha _{7}  \\
\alpha  _{4} & \alpha  _{5}  &   \alpha_{6} &  \alpha _{7} &  \alpha _{8}  \\
\alpha  _{5} & \alpha  _{6}  &   \alpha_{7} &  \alpha _{8} &  \alpha _{9}\\
\hline
\beta  _{1} & \beta  _{2} & \beta  _{3} & \beta_{4}  & \beta_{5}  
 \end{array}  \right) $$
is equal to \vspace{0.1 cm} \\

$$ 2^{i}\Gamma _{i}^{5 \times 5}+
 (2^{5}-2^{i-1})\cdot\Gamma _{i-1}^{5\times 5}\quad for\quad 0\leq i\leq 5, $$\vspace{0.1 cm}\\
\begin{equation*}
=  \begin{cases}
1 & \text{if  } i = 0, \\
37 &  \text{if  }    i=1, \\
138 & \text{if   } i = 2, \\
720 & \text{if   } i = 3, \\
4224  & \text{if   } i = 4. \\
11264  & \text{if   } i = 5. \\
\end{cases}
\end{equation*}
\vspace{0.1 cm}
\textbf{Application:}
 The number $ R_{q} $ of solutions
  $(Y_1,Z_1,U_{1}, Y_2,Z_2,U_{2},\ldots , Y_q,Z_q,U_{q} ) \in \big(\mathbb{F}_{2}[T] \big)^{3q}  $ \vspace{0.5 cm}\\
 of the polynomial equations\\
  \[\left\{\begin{array}{c}
 Y_{1}Z_{1} +Y_{2}Z_{2}+ \ldots + Y_{q}Z_{q} = 0  \\
   Y_{1}U_{1} + Y_{2}U_{2} + \ldots  + Y_{q}U_{q} = 0  
   \end{array}\right.\]
 satisfying the degree conditions \\
                   $$  degY_j \leq 4 , \quad degZ_j \leq 4 ,
                   \quad degU_{j}\leq 0 , \quad \text{for}  \quad 1\leq j \leq q $$ \\
  is connected with the numbers  $\Gamma _{i}^{\Big[\substack{5 \\ 1 }\Big] \times 5}$ 
   , in the way that the number $ R_{q} $ can be written as a linear combination of the $\Gamma _{i}^{\Big[\substack{5 \\ 1 }\Big] \times 5}$\\                  
  More precisely:\\
    \begin{align*}
& R_{q} 
= 2^{11q-14}\sum_{i = 0}^{5}
   \Gamma _{i}^{\left[\stackrel{5}{1}\right]\times 5}2^{-iq}\\ 
 &  = \left\{\begin{array}{ccc}
   95 & if & q=1,\\
   14752 & if & q=2, \\
   2^{11q-14}\left[1+ 37\cdot2^{-q} +138\cdot2^{-2q} + 720\cdot2^{-3q} +4224\cdot2^{-4q} +11264\cdot2^{-5q} \right]
 & if & 3\leq q.
  \end{array}\right. 
  \end{align*}
  \end{example}
\newpage  
  \begin{example}

 The number $ \Gamma_{i}^{\left[2\atop (2) \right]\times 4} $  of rank i matrices of the form \\

  \begin{displaymath}
 \left (  \begin{array} {cccc}
     \alpha  _{1} & \alpha  _{2}  &   \alpha_{3} &  \alpha _{4}  \\
\alpha  _{4} & \alpha  _{5}  &   \alpha_{6} &  \alpha _{7}   \\                                                                                             
\alpha  _{7} & \alpha  _{8}  &   \alpha_{9} &  \alpha _{10}  \\       
\alpha  _{10} & \alpha  _{11}  &   \alpha_{12} &  \alpha _{13}   
\end{array} \right ) \;    \overset{\text{rank}}{\sim}   \;
 \left (  \begin{array} {cccc}
\alpha  _{1} & \alpha  _{2}  &   \alpha_{3} &  \alpha _{4}   \\
\alpha  _{2} & \alpha  _{3}  &   \alpha_{4} &  \alpha _{5}  \\ 
\hline
\beta  _{11} & \beta  _{12} & \beta  _{13} & \beta_{14} \\
 \beta  _{21} & \beta  _{22} & \beta  _{23} & \beta_{24} 
\end{array} \right )    
\end{displaymath} \\
is equal to \vspace{0.1 cm} \\

 $2^{2i}\cdot\Gamma_{i}^{2\times 4}+3\cdot2^{i-1}\cdot(2^{4}-2^{i-1})\cdot \Gamma_{i-1}^{2\times 4} +(2^{4}-2^{i-1})\cdot(2^{4}-2^{i-2})\cdot\Gamma_{i-2}^{2\times 4}\quad \text{for} \;  0\leqslant i \leqslant 4$

\begin{equation*}
=  \begin{cases}
1 & \text{if  } i = 0, \\
57 &  \text{if  }    i=1, \\
910 & \text{if   } i = 2, \\
4536 & \text{if   } i = 3, \\
2688  & \text{if   } i = 4. \\
\end{cases}
\end{equation*}
\vspace{0.1 cm}
\textbf{Application:}
 The number $ R_{q} $ of solutions
  $(Y_1,Z_1,U_{1},V_{1}, Y_2,Z_2,U_{2},V_{2},\ldots , Y_q,Z_q,U_{q},V_{q} ) \in \big(\mathbb{F}_{2}[T] \big)^{4q}  $ \vspace{0.5 cm}\\
 of the polynomial equations\\
  \[\left\{\begin{array}{c}
 Y_{1}Z_{1} +Y_{2}Z_{2}+ \ldots + Y_{q}Z_{q} = 0  \\
   Y_{1}U_{1} + Y_{2}U_{2} + \ldots  + Y_{q}U_{q} = 0 \\
    Y_{1}V_{1} + Y_{2}V_{2} + \ldots  + Y_{q}V_{q} = 0 \\
     \end{array}\right.\]
 satisfying the degree conditions \\
                   $$  degY_j \leq 3 , \quad degZ_j \leq 1 ,
                   \quad degU_{j}\leq 0 ,  \quad degV_{j}\leq 0      \quad \text{for}  \quad 1\leq j \leq q $$ \\

 is equal to\\
 
  $   2^{8q-13}  \sum_{i=0}^{4}   \Gamma_{i}^{\left[2\atop (2) \right]\times 4}  \cdot2^{-iq} =
  2^{8q-13}\left[1+ 57\cdot2^{-q} +910\cdot2^{-2q} + 4536\cdot2^{-3q} +2688\cdot2^{-4q} \right] $

  \end{example}

  For more details see [2]

\section{Polynomial equations related to double persymmetric matrices}
\label{sec 3 }
  \begin{example}
  The most simple problem concerning double persymmetric matrices with entries in  $\mathbb{F}_{2}$:\\
  Compute the number   $ \Gamma_{i}^{\left[2\atop 2 \right]\times k} $ of double persymmetric matrices in  $\mathbb{F}_{2}$
  of rank  $0\leq i\leq \inf(4,k) $ of the form:
   $$   \left ( \begin{array} {cccccc}
 \alpha _{1} & \alpha _{2}  &  \ldots & \alpha _{k-1}  &  \alpha _{k} \\
\alpha _{2 } & \alpha _{3} &  \ldots  &  \alpha _{k} &  \alpha _{k+1} \\
\hline
\beta  _{1} & \beta  _{2}  & \ldots  &  \beta_{k-1} &  \beta _{k}  \\
\beta  _{2} & \beta  _{3}  & \ldots  &  \beta_{k} &  \beta _{k+1}  
 \end{array}  \right). $$ 
 We get 
 \begin{equation*}
 \Gamma_{i}^{\left[2\atop 2 \right]\times k} = 
\begin{cases}
1 & \text{if  } i = 0,\; k \geq 1,         \\
 9 & \text{if   } i=1,\; k > 1, \\
 3\cdot 2^{k+1} + 30 & \text{if   }  i = 2 ,\; k > 2, \\
  21\cdot 2^{k+1} -168 & \text{if   }  i=3 ,\; k > 3,\\ 
2^{2k+2} -3\cdot2^{k+4} + 128  & \text{if   }  i=4 ,\; k \geq 4
  \end{cases}    
   \end{equation*}
  \begin{equation*}
 \Gamma_{i}^{\left[2\atop 2 \right]\times i} = 
\begin{cases}
15 & \text{if  } i = 1,   \\
 54 & \text{if   } i=2, \\
168 & \text{if   }  i=3 ,\\ 
384 & \text{if   }  i=4 
  \end{cases}    
   \end{equation*}
\vspace{0.1 cm}
\textbf{Application:}
 The number $ R_{q} $ of solutions
  $(Y_1,Z_1,U_{1}, Y_2,Z_2,U_{2},\ldots , Y_q,Z_q,U_{q} ) \in \big(\mathbb{F}_{2}[T] \big)^{3q}  $ \vspace{0.5 cm}\\
 of the polynomial equations\\
  \[\left\{\begin{array}{c}
 Y_{1}Z_{1} +Y_{2}Z_{2}+ \ldots + Y_{q}Z_{q} = 0  \\
   Y_{1}U_{1} + Y_{2}U_{2} + \ldots  + Y_{q}U_{q} = 0  
   \end{array}\right.\]
 satisfying the degree conditions \\
                   $$  degY_j \leq k-1 , \quad degZ_j \leq 1 ,
                   \quad degU_{j}\leq 1 , \quad \text{for}  \quad 1\leq j \leq q $$ \\
  is connected with the numbers  $\Gamma _{i}^{\Big[\substack{2 \\ 2 }\Big] \times k}$ 
   , in the way that the number $ R_{q} $ can be written as a linear combination of the $\Gamma _{i}^{\Big[\substack{2 \\ 2 }\Big] \times k}$\\                  
  More precisely:\\
    \begin{align*}
& R_{q} 
= 2^{(k+4)q-2k-2}\sum_{i = 0}^{inf(4,k)}
   \Gamma _{i}^{\left[\stackrel{2}{2}\right]\times k}2^{-iq}\\ 
 &  = \left\{\begin{array}{ccc}
   2^{k} +15 & if & q=1,k\geqslant 4\\
2^{2k+8} +27\cdot2^{k+1}+192    & if & q=2,k\geqslant 4 \\
   2^{kq +4q -2k -2}\cdot \left[1+ 9\cdot2^{-q} +(3\cdot2^{k+1} +30)\cdot2^{-2q} + (21\cdot2^{k+1} -168)\cdot2^{-3q} \right] \\
   +  2^{kq +4q -2k -2}\cdot  \left[ (2^{2k+2} -3\cdot2^{k+4} +128)  \cdot 2^{-4q}  \right]  
 & if & 3\leq q, k\geqslant 4
  \end{array}\right. 
  \end{align*}
\textbf{The case k = 1}

\begin{equation*}
 R_{q} 
= 2^{5q-4}\sum_{i = 0}^{1}
   \Gamma _{i}^{\left[\stackrel{2}{2}\right]\times 1}2^{-iq} = 2^{5q-4} +15\cdot2^{4q-4}
\end{equation*}
\textbf{The case k = 2}

\begin{equation*}
 R_{q} 
= 2^{6q-6}\sum_{i = 0}^{2}
   \Gamma _{i}^{\left[\stackrel{2}{2}\right]\times 2}2^{-iq} = 2^{6q-6} +9\cdot2^{5q-6} + 54\cdot2^{4q-6}
\end{equation*}
\textbf{The case k = 3}

\begin{equation*}
 R_{q} 
= 2^{7q-8}\sum_{i = 0}^{3}
   \Gamma _{i}^{\left[\stackrel{2}{2}\right]\times 3}2^{-iq} = 2^{7q-8} +9\cdot2^{6q-8} + 78\cdot2^{5q-8} + 168\cdot2^{4q-8} 
\end{equation*}

  \end{example}

 \begin{example}
 The number  $ \Gamma_{i}^{\left[5\atop 5 \right]\times k} $ of double persymmetric  $10\times k$   matrices over  $\mathbb{F}_{2}$  of rank i of the form \\
 $$   \left ( \begin{array} {cccccc}
 \alpha _{1} & \alpha _{2}  &  \ldots & \alpha _{k-1}  &  \alpha _{k} \\
\alpha _{2 } & \alpha _{3} &  \ldots  &  \alpha _{k} &  \alpha _{k+1} \\
\alpha  _{3} & \alpha  _{4}  & \ldots  &  \alpha_{k+1} &  \alpha _{k+2}  \\
\alpha  _{4} & \alpha _{5}  & \ldots  &  \alpha_{k+2} &  \alpha _{k+3}  \\
\alpha  _{5} & \alpha  _{6}  & \ldots  &  \alpha_{k+3} &  \alpha _{k+4} \\ 
\hline
\beta  _{1} & \beta  _{2}  & \ldots  &  \beta_{k-1} &  \beta _{k}  \\
\beta  _{2} & \beta  _{3}  & \ldots  &  \beta_{k} &  \beta _{k+1}  \\
\beta  _{3} & \beta  _{4}  & \ldots  &  \beta_{k+1} &  \beta _{k+2}  \\
\beta  _{4} & \beta  _{5}  & \ldots  &  \beta_{k+2} &  \beta _{k+3}  \\
\beta  _{5} & \beta  _{6}  & \ldots  &  \beta_{k+3} &  \beta _{k+4}  
 \end{array}  \right). $$ 
 
is equal to \\
  \begin{equation*}
\begin{cases}
1 & \text{if  } i = 0,\; k \geq 1,         \\
 9 & \text{if   } i=1,\; k > 1, \\
  78 & \text{if   }  i = 2 ,\; k > 2, \\
 648 & \text{if   }  i=3 ,\; k > 3,\\ 
  5280 & \text{if   }  i=4 ,\; k > 4,\\
 3\cdot2^{k+4} + 39552 & \text{if   }  i=5 ,\; k > 5,\\ 
  21\cdot2^{k+4} + 290304& \text{if   }  i=6 ,\; k > 6,\\
   21\cdot2^{k+7} +1892352 & \text{if   }  i=7 ,\; k > 7,\\
    21\cdot2^{k+10} + 825753 6& \text{if   }  i=8 ,\; k > 8,\\
     21\cdot2^{k+13} - 44040192& \text{if   }  i=9 ,\; k > 9,\\
    2^{2k + 8}  -3\cdot 2^{k+16}  + 2^{25} & \text{if   }  i = 10 ,\; k  \geq 10.
 \end{cases}    
   \end{equation*}
   
   \textbf{The case k=i}
   
     \begin{equation*}
    \Gamma_{i}^{\left[5\atop 5 \right]\times i} =
     \begin{cases}
 2^{10} -1& \text{if   } i=1, \\
  2^{12} -10& \text{if   } i=2, \\
   2^{14} -88 & \text{if   } i=3, \\
    2^{16} -736 & \text{if   } i=4, \\
2^{18} - 6016 & \text{if   } i=5, \\
 2^{20} - 48640 & \text{if   } i=6, \\
  2^{22} - 385024& \text{if   } i=7, \\
      2^{24} - 3014656 & \text{if   } i=8, \\
   2^{26} - 23068672 & \text{if   } i=9, \\
     2^{28} - 167772160 & \text{if   } i=10, 
      \end{cases}    
   \end{equation*}

\textbf{The case k=4}

 The number  $ \Gamma_{i}^{\left[5\atop 5 \right]\times 4} $ of double persymmetric  $10\times 4$   matrices over  $\mathbb{F}_{2}$  of rank i of the form \\
 $$   \left ( \begin{array} {cccccc}
 \alpha _{1} & \alpha _{2}  &  \alpha _{3}  &  \alpha _{4} \\
\alpha _{2 } & \alpha _{3} &  \alpha _{4} &  \alpha _{5} \\
\alpha  _{3} & \alpha  _{4}   &  \alpha_{5} &  \alpha _{6}  \\
\alpha  _{4} & \alpha _{5}   &  \alpha_{6} &  \alpha _{7}  \\
\alpha  _{5} & \alpha  _{6}  &  \alpha_{7} &  \alpha _{8} \\ 
\hline
\beta  _{1} & \beta  _{2}  &  \beta_{3} &  \beta _{4}  \\
\beta  _{2} & \beta  _{3}   &  \beta_{4} &  \beta _{5}  \\
\beta  _{3} & \beta  _{4}  &  \beta_{5} &  \beta _{6}  \\
\beta  _{4} & \beta  _{5}   &  \beta_{6} &  \beta _{7}  \\
\beta  _{5} & \beta  _{6}   &  \beta_{7} &  \beta _{8}  
 \end{array}  \right). $$ 
 
is equal to \\
  \begin{equation*}
   \begin{cases}
1 & \text{if  } i = 0,  \\
 9 & \text{if   } i=1 \\
  78 & \text{if   }  i = 2 , \\
 648 & \text{if   }  i=3 \\ 
 64800  & \text{if   }  i=4 
 \end{cases}    
   \end{equation*}
   
  \textbf{Application:}
 The number $ R_{q} $ of solutions
  $(Y_1,Z_1,U_{1}, Y_2,Z_2,U_{2},\ldots , Y_q,Z_q,U_{q} ) \in \big(\mathbb{F}_{2}[T] \big)^{3q}  $ \vspace{0.5 cm}\\
 of the polynomial equations\\
  \[\left\{\begin{array}{c}
 Y_{1}Z_{1} +Y_{2}Z_{2}+ \ldots + Y_{q}Z_{q} = 0  \\
   Y_{1}U_{1} + Y_{2}U_{2} + \ldots  + Y_{q}U_{q} = 0  
   \end{array}\right.\]
 satisfying the degree conditions \\
                   $$  degY_j \leq 3 , \quad degZ_j \leq 4 ,
                   \quad degU_{j}\leq 4 , \quad \text{for}  \quad 1\leq j \leq q $$ \\
  is connected with the numbers  $\Gamma _{i}^{\Big[\substack{5 \\ 5 }\Big] \times 4}$ 
   , in the way that the number $ R_{q} $ can be written as a linear combination of the $\Gamma _{i}^{\Big[\substack{5 \\ 5 }\Big] \times 4}$\\          
  More precisely:\\
    \begin{align*}
& R_{q} 
= 2^{14q-16}\sum_{i = 0}^{4}
   \Gamma _{i}^{\left[\stackrel{5}{5}\right]\times 4}2^{-iq}\\ 
   & =  2^{14q-16}\cdot[1+9\cdot2^{-q} + 78\cdot2^{-2q} +648\cdot2^{-3q} + 64800\cdot2^{-4q} ]
  \end{align*}

 \end{example}   
\begin{example}
The fraction of square double persymmetric matrices which are invertible is equal to $ {3\over 8} $
\end{example}

\begin{example}  
  The number  $ \Gamma_{i}^{\left[2\atop 2+3 \right]\times 4} $ of double persymmetric  $7\times 4$   matrices over  $\mathbb{F}_{2}$  of rank i of the form \\
 $$   \left ( \begin{array} {cccc}
 \alpha _{1} & \alpha _{2} & \alpha _{3}  &  \alpha _{4} \\
\alpha _{2 } & \alpha _{3} &    \alpha _{4} &  \alpha _{5} \\
\hline
\beta  _{1} & \beta  _{2} &  \beta_{3} &  \beta _{4}  \\
\beta  _{2} & \beta  _{3}   &  \beta_{4} &  \beta _{5}  \\
\beta  _{3} & \beta  _{4}  &  \beta_{5} &  \beta _{6}  \\
\beta  _{4} & \beta  _{5}   &  \beta_{6} &  \beta _{7}  \\
\beta  _{5} & \beta  _{6}  &  \beta_{7} &  \beta _{8}  
 \end{array}  \right). $$ 
 is given by\\

   \begin{equation*}
\begin{cases}
1 & \text{if  } i = 0,        \\
 9 & \text{if   } i=1,\\
 94 & \text{if   }  i = 2,  \\
  600 & \text{if   }  i=3 ,\\
  7488 & \text{if   }  i=4 
   \end{cases}    
   \end{equation*}
   \textbf{Application:}\\
   
    The number $ R_{q} $ of solutions
  $(Y_1,Z_1,U_{1}, Y_2,Z_2,U_{2},\ldots , Y_q,Z_q,U_{q} ) \in \big(\mathbb{F}_{2}[T] \big)^{3q}  $ \vspace{0.5 cm}\\
 of the polynomial equations\\
  \[\left\{\begin{array}{c}
 Y_{1}Z_{1} +Y_{2}Z_{2}+ \ldots + Y_{q}Z_{q} = 0  \\
   Y_{1}U_{1} + Y_{2}U_{2} + \ldots  + Y_{q}U_{q} = 0  
   \end{array}\right.\]
 satisfying the degree conditions \\
                   $$  degY_j \leq 3 , \quad degZ_j \leq 1,
                   \quad degU_{j}\leq 4 , \quad \text{for}  \quad 1\leq j \leq q $$ \\
  is connected with the numbers  $\Gamma _{i}^{\Big[\substack{2 \\ 2+3 }\Big] \times 4}$ 
   , in the way that the number $ R_{q} $ can be written as a linear combination of the $\Gamma _{i}^{\Big[\substack{2 \\ 2+3 }\Big] \times 4}$\\          
  More precisely:\\
    \begin{align*}
& R_{q} 
= 2^{11q-13}\sum_{i = 0}^{4}
   \Gamma _{i}^{\left[\stackrel{2}{2+3}\right]\times 4}2^{-iq}\\ 
   & =  2^{11q-13}\cdot[1+9\cdot2^{-q} + 94\cdot2^{-2q} +600\cdot2^{-3q} + 7488\cdot2^{-4q} ]
  \end{align*}
   \end{example}

   \begin{example}  
  The number  $ \Gamma_{i}^{\left[2\atop {2+3\atop (1)} \right]\times 4} $ of rank i  $8\times 4$   matrices over  $\mathbb{F}_{2}$  of the form \\
 $$   \left ( \begin{array} {cccc}
 \alpha _{1} & \alpha _{2} & \alpha _{3}  &  \alpha _{4} \\
\alpha _{2 } & \alpha _{3} &    \alpha _{4} &  \alpha _{5} \\
\hline
\beta  _{1} & \beta  _{2} &  \beta_{3} &  \beta _{4}  \\
\beta  _{2} & \beta  _{3}   &  \beta_{4} &  \beta _{5}  \\
\beta  _{3} & \beta  _{4}  &  \beta_{5} &  \beta _{6}  \\
\beta  _{4} & \beta  _{5}   &  \beta_{6} &  \beta _{7}  \\
\beta  _{5} & \beta  _{6}  &  \beta_{7} &  \beta _{8} \\
\hline
 \gamma_{11} &  \gamma_{12} &  \gamma_{13} &  \gamma_{4} 
 \end{array}  \right). $$ 
 is equal to \\ 
 
$$ 2^{i}\cdot \Gamma _{i}^{\Big[\substack{2 \\ 2+3 }\Big] \times 4}+(2^{4} -2^{i-1})\cdot \Gamma _{i-1}^{\Big[\substack{2 \\ 2+3 }\Big] \times 4} 
\text{for} \quad  1\leqslant i \leqslant 4 $$

   \begin{equation*}
  =  \begin{cases}
1 & \text{if  } i = 0,        \\
 33 & \text{if   } i=1,\\
 502 & \text{if   }  i = 2,  \\
  5928 & \text{if   }  i=3 ,\\
  124608 & \text{if   }  i=4 
   \end{cases}    
   \end{equation*}
   
  \textbf{Application:}
  
     The number $ R_{q} $ of solutions
    $(Y_1,Z_1,U_{1},V_{1}, \ldots,Y_q,Z_q,U_{q},V_{q})\in \big( \mathbb{F}_{2}[T ]  \big)^{4q} $  
    of the polynomial equations\\ 
      \[\left\{\begin{array}{c}
 Y_{1}Z_{1} +Y_{2}Z_{2}+ \ldots + Y_{q}Z_{q} = 0,  \\
   Y_{1}U_{1} + Y_{2}U_{2} + \ldots  + Y_{q}U_{q} = 0,\\
    Y_{1}V_{1} + Y_{2}V_{2} + \ldots  + Y_{q}V_{q} = 0,\\   
 \end{array}\right.\]
  satisfying the degree conditions \\
                   $$  degY_j \leq 3 , \quad degZ_j \leq 1,
                   \quad degU_{j}\leq 4 ,  \quad degV_{j}\leq 0  \quad \text{for}  \quad 1\leq j \leq q $$ \\
  is connected with the numbers  $ \Gamma_{i}^{\left[2\atop {2+3\atop (1)} \right]\times 4} $
   , in the way that the number $ R_{q} $ can be written as a linear combination of the    $ \Gamma_{i}^{\left[2\atop {2+3\atop( 1)} \right]\times 4} $

   More precisely:\\
    \begin{align*}
& R_{q} 
= 2^{12q-17}\sum_{i = 0}^{4}
 \Gamma_{i}^{\left[2\atop {2+3\atop (1)} \right]\times 4} 2^{-iq}\\ 
  & =  2^{12q-17}\cdot[1+33\cdot2^{-q} + 502\cdot2^{-2q} +5928\cdot2^{-3q} + 124608\cdot2^{-4q} ]
  \end{align*}
 \end{example}

   \begin{example}

 The number $ \Gamma_{i}^{\left[2\atop {2\atop (4)} \right]\times 4} $ of rank i matrices of the form \\

  \begin{displaymath}
 \left (  \begin{array} {cccccccc}
     \alpha  _{1} & \alpha  _{2}  &   \alpha_{3} &  \alpha _{4}  &  \alpha  _{5} & \alpha  _{6}  &   \alpha_{7} &  \alpha _{8}  \\
        \alpha  _{7} & \alpha  _{8}  &   \alpha_{9} &  \alpha _{10}  &  \alpha  _{11} & \alpha  _{12}  &   \alpha_{13} &  \alpha _{14}  \\
           \alpha  _{13} & \alpha  _{14}  &   \alpha_{15} &  \alpha _{16}  &  \alpha  _{17} & \alpha  _{18}  &   \alpha_{19} &  \alpha _{20}  \\
       \alpha  _{19} & \alpha  _{20}  &   \alpha_{21} &  \alpha _{22}  &  \alpha  _{23} & \alpha  _{24}  &   \alpha_{25} &  \alpha _{26}       
           \end{array} \right ) \;   \overset{\text{rank}}{\sim}  \;
 \left (  \begin{array} {cccc}
\alpha  _{1} & \alpha  _{2}  &   \alpha_{3} &  \alpha _{4}   \\
\alpha  _{2} & \alpha  _{3}  &   \alpha_{4} &  \alpha _{5}  \\ 
\hline \\
\beta  _{1} & \beta _{2}  &   \beta_{3} &  \beta _{4}   \\
\beta  _{2} & \beta  _{3}  &   \beta_{4} &  \beta _{5}  \\ 
\hline\\
\gamma _{11} & \gamma _{12} & \gamma  _{13} & \gamma_{14} \\
 \gamma _{21} & \gamma  _{22} & \gamma _{23} & \gamma_{24} \\
 \gamma _{31} & \gamma _{32} & \gamma  _{33} & \gamma_{34} \\
 \gamma _{41} & \gamma  _{42} & \gamma _{43} & \gamma_{44} 
\end{array} \right )    
\end{displaymath} \\
is equal to \vspace{0.1 cm} \\

    $ 2^{4i}\cdot\Gamma _{i}^{\Big[\substack{2 \\ 2 }\Big] \times 4}
+15\cdot2^{(i-1)3}(2^{k}-2^{i-1})\cdot\Gamma _{i-1}^{\Big[\substack{2 \\ 2 }\Big] \times 4}\\
  +   35\cdot2^{2i-4}(2^{k}-2^{i-1})(2^{k}-2^{i-2})\cdot \Gamma _{i-2}^{\Big[\substack{2 \\ 2 }\Big] \times 4}   \\
      + 15\cdot2^{i-3}(2^{k}-2^{i-1})(2^{k}-2^{i-2})(2^{k}-2^{i-3})\Gamma _{i-3}^{\Big[\substack{2 \\ 2 }\Big] \times 4} \\
    + (2^{k}-2^{i-1})(2^{k}-2^{i-2})(2^{k}-2^{i-3})(2^{k}-2^{i-4})\Gamma _{i-4}^{\Big[\substack{2 \\ 2 }\Big] \times 4} 
        for\quad 0\leq i\leq inf(4,8), $\vspace{0.5 cm} \\

\begin{equation*}
=  \begin{cases}
1 & \text{if  } i = 0, \\
369 &  \text{if  }    i=1, \\
54726 & \text{if   } i = 2, \\
3765384 & \text{if   } i = 3, \\
63288384 & \text{if   } i = 4. \\
\end{cases}
\end{equation*}
\vspace{0.1 cm}

Then the number of solutions \\
 $(Y_1,Z_1,U_{1},V_{1}^{(1)},V_{2}^{(1)}, V_{3}^{(1)},V_{4}^{(1)}, Y_2,Z_2,U_{2},V_{1}^{(2)},V_{2}^{(2)}, 
V_{3}^{(2)},V_{4}^{(2)},\ldots  Y_q,Z_q,U_{q},V_{1}^{(q)},V_{2}^{(q)}, V_{3}^{(q)},V_{4}^{(q)}   )\in  \big( \mathbb{F}_{2}[T ]  \big)^{7q}  $ \vspace{0.5 cm}\\

 of the polynomial equations  \vspace{0.5 cm}
  \[\left\{\begin{array}{c}
 Y_{1}Z_{1} +Y_{2}Z_{2}+ \ldots + Y_{q}Z_{q} = 0  \\
  Y_{1}U_{1} +Y_{2}U_{2}+ \ldots + Y_{q}U_{q} = 0 \\
   Y_{1}V_{1}^{(1)} + Y_{2}V_{1}^{(2)} + \ldots  + Y_{q}V_{1}^{(q)} = 0  \\
    Y_{1}V_{2}^{(1)} + Y_{2}V_{2}^{(2)} + \ldots  + Y_{q}V_{2}^{(q)} = 0\\
      Y_{1}V_{3}^{(1)} + Y_{2}V_{3}^{(2)} + \ldots  + Y_{q}V_{3}^{(q)} = 0\\
      Y_{1}V_{4}^{(1)} + Y_{2}V_{4}^{(2)} + \ldots  + Y_{q}V_{4}^{(q)} = 0 
 \end{array}\right.\]
    satisfying the degree conditions \\
                   $$  \deg Y_i \leq 3 , \quad \deg Z_i \leq 1 ,\quad \deg U_i \leq 1 ,
                   \quad \deg V_{j}^{i} \leq 0 , \quad  for \quad 1\leq j\leq 4  \quad 1\leq i \leq q $$ \\
                   
  is equal to 
$$  2^{12q - 26}\sum_{i = 0}^{4}
  \Gamma_{i}^{\left[ 2\atop{ 2\atop (4)} \right]\times k} 2^{-iq}= 2^{12q-26}\left[1+ 369\cdot2^{-q} +54726\cdot2^{-2q} + 3765384\cdot2^{-3q} +
63288384 \cdot2^{-4q} \right] $$

  \end{example}

 For more details see [3], [4]

   \newpage
\section{Polynomial equations related to triple  persymmetric matrices}
\label{sec 4}
\begin{example}

  The number  $ \Gamma_{i}^{\left[1\atop {1\atop 1} \right]\times 1} $ of triple  persymmetric  $3\times 1$   matrices over  $\mathbb{F}_{2}$  of rank i of the form \\
 $$   \left ( \begin{array} {cccccc}
\alpha _{1}   \\
\hline \\
\beta  _{1}    \\  
\hline \\
\gamma  _{1}   \\
 \end{array}  \right). $$ 
 is equal to
 \[ \begin{cases}
1  &\text{if  }  i = 0 \\
7  & \text{if   } i = 2 \\
\end{cases}
\]

\textbf{Application:}

The number of solutions    $(y_1,z_1,u_{1},v_{1}, \ldots,y_q,z_q,u_{q},v_{q})  \in    \mathbb{F}_{2}^{4q}  $ of the following system 
of quadratic equations : \\
   \[\left\{\begin{array}{c}
 y_{1}z_{1} +y_{2}z_{2}+ \ldots + y_{q}z_{q} = 0,  \\
   y_{1}u_{1} + y_{2}u_{2} + \ldots  + y_{q}u_{q} = 0,\\
    y_{1}v_{1} + y_{2}v_{2} + \ldots  + y_{q}v_{q} = 0,\\  
 \end{array}\right.\]
    is equal to 
\begin{align*}
&  2^{4q-3}\cdot \sum_{i = 0}^{1} \Gamma_{i}^{\left[1\atop{ 1\atop 1} \right]\times 1}  \cdot2^{- qi} \\
& = 2^{4q-3}\cdot\big(1 + 7\cdot2^{-q} + \big ) \\
 & =   2^{4q-3} + 7\cdot2^{3q-3} 
\end{align*}

\textbf{Generalization} 

  The number  $ \Gamma_{i}^{\left[1\atop{\vdots \atop 1}\right]\times 3} $  of n-times   persymmetric  $n\times 1$   matrices over  $\mathbb{F}_{2}$  of rank i of the form \\
 $$   \left ( \begin{array} {cccccc}
\alpha _{1}   \\
\hline \\
\alpha _{2}    \\  
\hline \\
\vdots \\
\hline \\
\alpha  _{n}   \\
 \end{array}  \right). $$ 
 is equal to
 \[ \begin{cases}
1  &\text{if  }  i = 0 \\
2^n -1  & \text{if   } i = 2 \\
\end{cases}
\]
\textbf{Application:}

The number of solutions    $(y_1,z_1^{(1)},z_{1}^{(2)}, \ldots, z_{1}^{(n)}, \ldots,  y_q,z_{q}^{(1)},z_{q}^{(2)}, \ldots, z_{q}^{(n)}  )  \in    \mathbb{F}_{2}^{(n+1)q}  $ of the following system 
of quadratic equations : \\
   \[\left\{\begin{array}{c}
 y_{1}z_{1}^{(1)} +y_{2}z_{2}^{(1)}+ \ldots + y_{q}z_{q}^{(1)} = 0,  \\
  y_{1}z_{1}^{(2)} +y_{2}z_{2}^{(2)}+ \ldots + y_{q}z_{q}^{(2)} = 0,  \\
   y_{1}z_{1}^{(3)} +y_{2}z_{2}^{(3)}+ \ldots + y_{q}z_{q}^{(3)} = 0,  \\
   \vdots \\
 y_{1}z_{1}^{(n)} +y_{2}z_{2}^{(n)}+ \ldots + y_{q}z_{q}^{(n)} = 0
 \end{array}\right.\]
    is equal to 
\begin{align*}
&  2^{(n+1)q-n}\cdot \sum_{i = 0}^{1}  \Gamma_{i}^{\left[1\atop{\vdots \atop 1}\right]\times 3}\cdot2^{- qi} \\
& = 2^{(n+1)q-n}\cdot\big(1 + (2^n -1)\cdot2^{-q}\big )  = 2^{nq-n}\cdot[2^q+2^n-1]
\end{align*}

\end{example}

  \begin{example}  
  The number  $ \Gamma_{i}^{\left[2\atop {2\atop 2} \right]\times 6} $ of triple  persymmetric  $6\times 6$   matrices over  $\mathbb{F}_{2}$  of rank i of the form \\
 $$   \left ( \begin{array} {cccccc}
\alpha _{1} & \alpha _{2}  & \alpha _{3} & \alpha _{4}  & \alpha _{5}  &  \alpha _{6} \\
\alpha _{2 } & \alpha _{3} &  \alpha _{4} & \alpha _{5} &  \alpha _{6} &  \alpha _{7} \\
\hline
\beta  _{1} & \beta  _{2}  & \beta  _{3}& \beta  _{4}   &  \beta_{5} &  \beta _{6}  \\
\beta  _{2} & \beta  _{3}  &  \beta  _{3}&  \beta  _{5} &  \beta_{6} &  \beta _{7}  \\
\hline
\gamma  _{1} & \gamma   _{2}  & \gamma   _{3}& \gamma   _{4}  & \gamma  _{5} &  \gamma  _{6}  \\
\gamma  _{2} & \gamma  _{3} & \gamma   _{4}  & \gamma   _{5}  & \gamma  _{6} &  \gamma  _{7}  \\
 \end{array}  \right). $$ 
 is equal to
 \[ \begin{cases}
1  &\text{if  }  i = 0 \\
21   &\text{if  }  i = 1 \\
1162     & \text{if   } i = 2 \\
20160          & \text{if   } i = 3 \\
258720           & \text{if   } i = 4 \\
1128960           & \text{if   } i = 5 \\
688128  & \text{if  } i = 6
\end{cases}
\]
\textbf{Application:}

  The number $ R_{q} $ of solutions \\
 $(Y_1,Z_1,U_{1},V_{1}, \ldots,Y_q,Z_q,U_{q},V_{q})  \in \big( \mathbb{F}_{2}[T ]  \big)^{4q} $ of the polynomial equations
   \[\left\{\begin{array}{c}
 Y_{1}Z_{1} +Y_{2}Z_{2}+ \ldots + Y_{q}Z_{q} = 0,  \\
   Y_{1}U_{1} + Y_{2}U_{2} + \ldots  + Y_{q}U_{q} = 0,\\
    Y_{1}V_{1} + Y_{2}V_{2} + \ldots  + Y_{q}V_{q} = 0,\\  
 \end{array}\right.\]
  satisfying the degree conditions \\
                   $$  degY_j \leq 5 , \quad degZ_j \leq 1 ,\quad degU_{j}\leq 1,\quad degV_{j}\leq 1 \quad for \quad 1\leq j \leq q. $$ \\                           
is equal to 
\begin{align*}
&  2^{12q-21}\cdot \sum_{i = 0}^{6} \Gamma_{i}^{\left[2\atop{ 2\atop 2} \right]\times 6}  \cdot2^{- qi} \\
& = 2^{12q-21}\cdot\big(1 + 21\cdot2^{-q} + 1162\cdot2^{-2q} + 20160\cdot2^{-3q} + 258720\cdot2^{-4q} + 1128960 \cdot2^{-5q}
 + 688128 \cdot2^{-6q}\big ) \\
 & =   2^{6q-21}\cdot\big(2^{6q} + 21\cdot2^{5q} + 1162\cdot2^{4q} + 20160\cdot2^{3q} + 258720\cdot2^{2q} + 1128960 \cdot2^{q}
 + 688128 \big ) 
\end{align*}
\end{example}
\begin{example}
The fraction of square triple persymmetric matrices which are invertible is equal to $ {21\over 64} $
\end{example}

\begin{example}

  The number  $ \Gamma_{i}^{\left[2\atop {2\atop 2} \right]\times 2} $ of triple  persymmetric  $6\times 2$   matrices over  $\mathbb{F}_{2}$  of rank i of the form \\
 $$   \left ( \begin{array} {cccccc}
\alpha _{1} & \alpha _{2}   \\
\alpha _{2 } & \alpha _{3}  \\
\hline
\beta  _{1} & \beta  _{2}    \\
\beta  _{2} & \beta  _{3}   \\
\hline
\gamma  _{1} & \gamma   _{2}    \\
\gamma  _{2} & \gamma  _{3}  \\
 \end{array}  \right). $$ 
 is equal to
 \[ \begin{cases}
1  &\text{if  }  i = 0 \\
21   &\text{if  }  i = 1 \\
490    & \text{if   } i = 2 \\
\end{cases}
\]

\textbf{Application:}

  The number $ R_{q} $ of solutions \\
 $(Y_1,Z_1,U_{1},V_{1}, \ldots,Y_q,Z_q,U_{q},V_{q})  \in \big( \mathbb{F}_{2}[T ]  \big)^{4q} $ of the polynomial equations
   \[\left\{\begin{array}{c}
 Y_{1}Z_{1} +Y_{2}Z_{2}+ \ldots + Y_{q}Z_{q} = 0,  \\
   Y_{1}U_{1} + Y_{2}U_{2} + \ldots  + Y_{q}U_{q} = 0,\\
    Y_{1}V_{1} + Y_{2}V_{2} + \ldots  + Y_{q}V_{q} = 0,\\  
 \end{array}\right.\]
  satisfying the degree conditions \\
                   $$  degY_j \leq 1 , \quad degZ_j \leq 1 ,\quad degU_{j}\leq 1,\quad degV_{j}\leq 1 \quad for \quad 1\leq j \leq q. $$ \\                           
is equal to 
\begin{align*}
&  2^{8q-9}\cdot \sum_{i = 0}^{2} \Gamma_{i}^{\left[2\atop{ 2\atop 2} \right]\times 2}  \cdot2^{- qi} \\
& = 2^{8q-9}\cdot\big(1 + 21\cdot2^{-q} + 490\cdot2^{-2q} \big ) \\
 & =   2^{8q-9} + 21\cdot2^{7q-9} + 490\cdot2^{6q-9} 
\end{align*}
\textbf{Application :}
The number of solutions $ (x_1,x_2,\ldots, x_{15},x_{16}) \in  \mathbb{F}_{2}^{16} $ of the following system 
of quadratic equations : \\

  \[\left\{\begin{array}{c}
 x_{1}x_{5} +x_{3}x_{7} = 0,  \\
   x_{1}x_{9} +x_{3}x_{11}   = 0,\\
   x_{1}x_{13} +x_{3}x_{15} = 0,\\ 
 x_{1}x_{6} +x_{2}x_{5} + x_{3}x_{8} +x_{4}x_{7} = 0\\  
  x_{1}x_{10} +x_{2}x_{9} + x_{3}x_{12} +x_{4}x_{11} = 0\\ 
     x_{1}x_{14} +x_{2}x_{13} + x_{3}x_{16} +x_{4}x_{15} = 0\\ 
      x_{2}x_{6} +x_{4}x_{8} = 0,  \\
   x_{2}x_{10} +x_{4}x_{12}   = 0,\\
   x_{2}x_{14} +x_{4}x_{16} = 0,\\   
    \end{array}\right.\]
    
   is equal to $R_{2} =   2^{7} + 21\cdot2^{5} + 490\cdot2^{3} = 4720 $\\
  \begin{proof} 
  Set 
   \[\left\{\begin{array}{c}
   Y_{1} = x_{1} +x_{2} \cdot T \\
     Y_{2} = x_{3} +x_{4} \cdot T \\ 
       Z_{1} = x_{5} +x_{6} \cdot T \\
     Z_{2} = x_{7} +x_{8} \cdot T \\ 
        U_{1} = x_{9} +x_{10} \cdot T \\
     U_{2} = x_{11} +x_{12} \cdot T \\ 
         V_{1} = x_{13} +x_{14} \cdot T \\
     V_{2} = x_{15} +x_{16} \cdot T 
     \end{array}\right.\]

  Then we obtain 
 
     \begin{equation*}
   \begin{cases}
  Y_{1}Z_{1} +Y_{2}Z_{2}= 0,  \\
   Y_{1}U_{1} + Y_{2}U_{2} = 0,\\    
       Y_{1}V_{1} + Y_{2}V_{2}  = 0 
  \end{cases} \; \Leftrightarrow  \; \begin{cases}
    x_{1}x_{5} +x_{3}x_{7} = 0,  \\
   x_{1}x_{9} +x_{3}x_{11}   = 0,\\
   x_{1}x_{13} +x_{3}x_{15} = 0,\\ 
 x_{1}x_{6} +x_{2}x_{5} + x_{3}x_{8} +x_{4}x_{7} = 0\\  
  x_{1}x_{10} +x_{2}x_{9} + x_{3}x_{12} +x_{4}x_{11} = 0\\ 
     x_{1}x_{14} +x_{2}x_{13} + x_{3}x_{16} +x_{4}x_{15} = 0\\ 
      x_{2}x_{6} +x_{4}x_{8} = 0,  \\
   x_{2}x_{10} +x_{4}x_{12}   = 0,\\
   x_{2}x_{14} +x_{4}x_{16} = 0,\\   
 \end{cases}
\end{equation*}

\end{proof}
 
 \end{example}
   \begin{example}

 The number $ \Gamma_{i}^{\left[2\atop {2\atop {2 \atop (3)}} \right]\times 4} $ of rank i matrices of the form \\

  \begin{displaymath}
 \left (  \begin{array} {ccccccccc}
     \alpha  _{1} & \alpha  _{2}  &   \alpha_{3} &  \alpha _{4}  &  \alpha  _{5} & \alpha  _{6}  &   \alpha_{7} &  \alpha _{8}  &  \alpha _{9}  \\
        \alpha  _{7} & \alpha  _{8}  &   \alpha_{9} &  \alpha _{10}  &  \alpha  _{11} & \alpha  _{12}  &   \alpha_{13} &  \alpha _{14} &  \alpha _{15}   \\
           \alpha  _{13} & \alpha  _{14}  &   \alpha_{15} &  \alpha _{16}  &  \alpha  _{17} & \alpha  _{18}  &   \alpha_{19} &  \alpha _{20} &  \alpha _{21}    \\
       \alpha  _{19} & \alpha  _{20}  &   \alpha_{21} &  \alpha _{22}  &  \alpha  _{23} & \alpha  _{24}  &   \alpha_{25} &  \alpha _{26}  &  \alpha _{27}       
           \end{array} \right ) \; \overset{\text{rank}}{\sim}  \;
 \left (  \begin{array} {cccc}
\alpha  _{1} & \alpha  _{2}  &   \alpha_{3} &  \alpha _{4}   \\
\alpha  _{2} & \alpha  _{3}  &   \alpha_{4} &  \alpha _{5}  \\ 
\hline \\
\beta  _{1} & \beta _{2}  &   \beta_{3} &  \beta _{4}   \\
\beta  _{2} & \beta  _{3}  &   \beta_{4} &  \beta _{5}  \\ 
\hline \\
\gamma _{1} & \gamma _{2} & \gamma  _{3} & \gamma_{4} \\
 \gamma _{2} & \gamma  _{3} & \gamma _{4} & \gamma_{5} \\
\hline \\
\delta_{11} & \delta_{12} & \delta  _{13} & \delta_{14} \\
\delta_{21} & \delta_{22} & \delta  _{23} & \delta_{24} \\
\delta_{31} & \delta_{32} & \delta  _{33} & \delta_{34} 
\end{array} \right )    
\end{displaymath} \\
is equal to \vspace{0.1 cm} \\

  $ 2^{3i}\cdot \Gamma_{i}^{\left[2\atop {2\atop 2} \right]\times 4}\\
    +7\cdot2^{(i-1)2}(2^{4}-2^{i-1})\cdot \Gamma_{i-1}^{\left[2\atop {2\atop 2} \right]\times 4}\\
  +   7\cdot2^{i-2}(2^{4}-2^{i-1})(2^{4}-2^{i-2})\cdot \Gamma_{i-2}^{\left[2\atop {2\atop 2} \right]\times 4}  \\
   + (2^{4}-2^{i-1})(2^{4}-2^{i-2})(2^{4}-2^{i-3})\cdot \Gamma_{i-3}^{\left[2\atop {2\atop 2} \right]\times 4}        for\quad 0\leq i\leq inf(4,9), $\vspace{0.5 cm} \\
Hence :\\

\begin{equation*}
 \Gamma_{i}^{\left[2\atop {2\atop {2 \atop (3)}} \right]\times 4} =  \begin{cases}
1 & \text{if  } i = 0, \\
273 &  \text{if  }    i=1, \\
41062 & \text{if   } i = 2, \\
3807048 & \text{if   } i = 3, \\
130369344 & \text{if   } i = 4. \\
\end{cases}
\end{equation*}
\vspace{0.1 cm}
\end{example}

For more details see [4]

 \section{The inverse problem}
\label{sec 5}

How to compute the number $ R_{q} $ of solutions \\
 $(Y_1,Z_1,U_{1},V_{1}, \ldots,Y_q,Z_q,U_{q},V_{q})  \in \big( \mathbb{F}_{2}[T ]  \big)^{4q} $ of the polynomial equations
   \[\left\{\begin{array}{c}
 Y_{1}Z_{1} +Y_{2}Z_{2}+ \ldots + Y_{q}Z_{q} = 0,  \\
   Y_{1}U_{1} + Y_{2}U_{2} + \ldots  + Y_{q}U_{q} = 0,\\
    Y_{1}V_{1} + Y_{2}V_{2} + \ldots  + Y_{q}V_{q} = 0,\\  
 \end{array}\right.\]
  satisfying the degree conditions \\
                   $$  degY_j \leq k-1 , \quad degZ_j \leq s-1 ,\quad degU_{j}\leq s+m-1,\quad degV_{j}\leq s+m+l-1 \quad for \quad 1\leq j\leq q. ?$$ \\

\textbf{Response: } \vspace{0,1 cm} \\
  We need only to compute                  
 the number  $\Gamma_{i}^{\left[s\atop{ s+m\atop s+m+l} \right]\times k} $ of triple persymmetric  $(3s+2m+l)\times k $ 
 rank i matrices 
 over $\mathbb{F}_{2}$ for $  0\leqslant i \leqslant \inf (k, 3s+2m+l) $ of the form
 $$   \left ( \begin{array} {cccccc}
\alpha _{1} & \alpha _{2}  &  \ldots & \alpha _{k-1}  &  \alpha _{k} \\
\alpha _{2 } & \alpha _{3} &  \ldots  &  \alpha _{k} &  \alpha _{k+1} \\
\vdots & \vdots & \vdots    & \vdots  &  \vdots \\
\alpha _{s-1} & \alpha _{s} & \ldots  &  \alpha _{s+k-3} &  \alpha _{s+k-2}  \\
 \alpha  _{s } & \alpha  _{s +1} & \ldots & \alpha  _{s +k-2}& \alpha  _{s +k-1}\\
\hline \\
\beta  _{1} & \beta  _{2}  & \ldots  &  \beta_{k-1} &  \beta _{k}  \\
\beta  _{2} & \beta  _{3}  & \ldots  &  \beta_{k} &  \beta _{k+1}  \\
\vdots & \vdots    &  \vdots & \vdots  &  \vdots \\
\beta  _{m+1} & \beta  _{m+2}  & \ldots  &  \beta_{k+m-1} &  \beta _{k+m}  \\
\vdots & \vdots    &  \vdots & \vdots  &  \vdots \\
\beta  _{s+m-1} & \beta  _{s+m}  & \ldots  &  \beta_{s+m+k-3} &  \beta _{s+m+k-2}  \\
 \beta _{s+m} & \beta _{s+m+1} & \ldots & \beta _{s+m+k-2} & \beta _{s+m+k-1}\\
\hline \\
\gamma  _{1} & \gamma   _{2}  & \ldots  & \gamma  _{k-1} &  \gamma  _{k}  \\
\gamma  _{2} & \gamma  _{3}  & \ldots  & \gamma  _{k} &  \gamma  _{k+1}  \\
\vdots & \vdots    &  \vdots & \vdots  &  \vdots \\
 \gamma  _{m+1} &  \gamma _{m+2}  & \ldots  & \gamma _{k+m-1} &  \gamma  _{k+m}  \\
\vdots & \vdots   &  \vdots & \vdots  &  \vdots \\
 \gamma  _{s+m-1} & \gamma  _{s+m}  & \ldots  & \gamma  _{s+m+k-3} &  \gamma  _{s+m+k-2}  \\
  \gamma  _{s+m} & \gamma  _{s+m+1}  & \ldots  & \gamma  _{s+m+k-2} &  \gamma  _{s+m+k-1}\\
\gamma  _{s+m+1} & \gamma  _{s+m+2}  & \ldots  & \gamma  _{s+m+k-1} &  \gamma  _{s+m+k}\\
\vdots & \vdots   &  \vdots & \vdots  &  \vdots \\
\gamma  _{s+m+l} & \gamma  _{s+m+l+1}  & \ldots  & \gamma  _{s+m+l+k-2} &  \gamma  _{s+m+l+k-1}
\end{array}  \right). $$ 
To compute those numbers we use the following reductions formulas:\\

    \begin{align*}
   \Gamma  _{2s+1+m+j}^{\left[s\atop{ s+m\atop s+m+l} \right]\times k} & = 16^{j}\cdot \Gamma  _{2s+1+m}^{\left[s\atop{ s+m\atop s+m+l-j} \right]\times (k-j)} 
   &  \text{if  } 0 \leq j \leq l,\;   k \geqslant 2s+1+m+j \\
      \Gamma  _{2s+1+m+l+j}^{\left[s\atop{ s+m\atop s+m+l} \right]\times k} & = 16^{2j+l}\cdot \Gamma  _{2s+1+m-j}^{\left[s\atop{ s+m-j\atop s+m-j} \right]\times (k-2j-l)} 
   &  \text{if  } 0 \leq j \leq m,  \; k \geqslant 2s+1+m+l+j     \\
\Gamma  _{2s+1+2m+l+j}^{\left[s\atop{ s+m\atop s+m+l} \right]\times k} & = 16^{2m+l+3j}\cdot \Gamma  _{2(s-j)+1}^{\left[s-j\atop{ s-j\atop s-j} \right]\times (k-2m-l-3j)} 
   &  \text{if  } 0 \leq j \leq s-1, \;  k \geqslant 2s+1+2m+l+j
\end{align*} 

Then $R_{q}$ is equal to a linear combination of the  $\Gamma_{i}^{\left[s\atop{ s+m\atop s+m+l} \right]\times k} \; \text{for} \;  0\leqslant i\leqslant \inf(k, 3s+2m+l)$ \\
  More precisely:\\
    \begin{align*}
& R_{q} 
= 2^{(k+3s+2m+l)q -(3k+3s+2m+l-3)} \sum_{i = 0}^{\inf(k, 3s+2m+l)}
 \Gamma_{i}^{\left[s\atop {s+m\atop s+m+l} \right]\times k} 2^{-iq}
   \end{align*}

\textbf{Remark : } \vspace{0.1 cm} \\
We have computed the  $\Gamma_{i}^{\left[s\atop{ s+m\atop s+m+l} \right]\times k}$ for  $ 0\leqslant i\leqslant \inf(k, 3s+2m+l)$ 
   in the case $l=0.$ [see (4)]
 \subsection{The inverse problem in the case $s=2,\;m =3,\;l\geqslant 4$ }
\label{subsec 1}
The results in this subsection are new.

 \begin{example}
To   compute the number $ R_{q} $ of solutions \\
 $(Y_1,Z_1,U_{1},V_{1}, \ldots,Y_q,Z_q,U_{q},V_{q})  \in \big( \mathbb{F}_{2}[T ]  \big)^{4q} $ of the polynomial equations
   \[\left\{\begin{array}{c}
 Y_{1}Z_{1} +Y_{2}Z_{2}+ \ldots + Y_{q}Z_{q} = 0,  \\
   Y_{1}U_{1} + Y_{2}U_{2} + \ldots  + Y_{q}U_{q} = 0,\\
    Y_{1}V_{1} + Y_{2}V_{2} + \ldots  + Y_{q}V_{q} = 0,\\  
 \end{array}\right.\]
  satisfying the degree conditions \\
                   $$  degY_j \leq k-1 , \quad degZ_j \leq 1 ,\quad degU_{j}\leq 4,\quad degV_{j}\leq 4+l \quad for \quad 1\leq j \leq q, $$ \\    
   we need only to compute                  
 the number  $\Gamma_{i}^{\left[2\atop{ 2+3\atop 2+3+l} \right]\times k} $ of triple persymmetric  $(12+l)\times k $ rank i  matrices 
 over $\mathbb{F}_{2}$ for $ 0\leqslant i\leqslant \inf(k, 12+l) $ of the form                 
 $$   \left ( \begin{array} {cccccc}
\alpha _{1} & \alpha _{2}  &  \ldots & \alpha _{k-1}  &  \alpha _{k} \\
\alpha _{2 } & \alpha _{3} &  \ldots  &  \alpha _{k} &  \alpha _{k+1} \\
\hline
\beta  _{1} & \beta  _{2}  & \ldots  &  \beta_{k-1} &  \beta _{k}  \\
\beta  _{2} & \beta  _{3}  & \ldots  &  \beta_{k} &  \beta _{k+1}  \\
\beta  _{3} & \beta  _{4}  & \ldots  &  \beta_{k+1} &  \beta _{k+2}  \\
\beta  _{4} & \beta  _{5}  & \ldots  &  \beta_{k+2} &  \beta _{k+3}  \\
\beta  _{5} & \beta  _{6}  & \ldots  &  \beta_{k+3} &  \beta _{k+4}  \\
\hline
\gamma  _{1} & \gamma   _{2}  & \ldots  & \gamma  _{k-1} &  \gamma  _{k}  \\
\gamma  _{2} & \gamma  _{3}  & \ldots  & \gamma  _{k} &  \gamma  _{k+1}  \\
\gamma  _{3} & \gamma  _{4}  & \ldots  & \gamma  _{k+1} &  \gamma  _{k+2}  \\
\vdots & \vdots    &  \vdots & \vdots  &  \vdots \\
 \gamma  _{5+l} & \gamma  _{6+l}  & \ldots  & \gamma  _{k+l+3} &  \gamma  _{k+l+4}  \\
 \end{array}  \right). $$ 
Proceeding as in (4), using the following reductions formulas:\\
   \begin{align*}
   \Gamma  _{8+j}^{\left[2\atop{ 2+3\atop 2+3+l} \right]\times k} & = 16^{j}\cdot \Gamma  _{8}^{\left[2\atop{ 2+3\atop 2+3+l-j} \right]\times (k-j)} 
   &  \text{if  } 0 \leq j \leq l, \; k \geqslant 8+j \\
      \Gamma  _{8+l+j}^{\left[2\atop{ 2+3\atop 2+3+l} \right]\times k} & = 16^{2j+l}\cdot \Gamma  _{8-j}^{\left[2\atop{ 2+3-j\atop 2+3-j} \right]\times (k-2j-l)} 
   &  \text{if  } 0 \leq j \leq 3, \; k \geqslant 8+l+j \\
\Gamma  _{11+l+j}^{\left[2\atop{ 2+3\atop 2+3+l} \right]\times k} & = 16^{6+l+3j}\cdot \Gamma  _{2(2-j)+1}^{\left[2-j\atop{ 2-j\atop 2-j} \right]\times (k-6-l-3j)} 
   &  \text{if  } 0 \leq j \leq 1, \; k \geqslant 11+l+j 
\end{align*} 
 we obtain:\\
  \begin{equation*}
\Gamma_{i}^{\left[2\atop{ 2+3\atop 2+3+l} \right]\times k} 
 = \begin{cases}
 1  & \text{if  } i = 0, \; k > 0\\
 21  & \text{if  } i = 1,  \; k > 1\\
 2^{k+1} + r_{2} (=362 ) & \text{if  } i = 2, \; k > 2 \\
 9\cdot2^{k+1} +  r_{3} (=6048) & \text{if  } i = 3, \; k > 3 \\
 39\cdot2^{k+2} + r_{4} (=98784) & \text{if  } i = 4,  \; k > 4\\
 97\cdot2^{k+4} +  r_{5}(= 1580288) & \text{if  } i = 5, \; k > 5 \\
 225\cdot2^{k+6} +  r_{6}(= 25135104)  & \text{if  } i = 6, \; k > 6 \\
 2^{2k+5} +417\cdot2^{k+8} +r_{7}(= 402571264=24571\cdot2^{14}   ) & \text{if   } i = 7,  \; k > 7
 \end{cases}
\end{equation*}

  \begin{equation*}
\Gamma_{8+i}^{\left[2\atop{ 2+3\atop 2+3+l} \right]\times k} 
 = \begin{cases}
 3\cdot2^{2k+5+2i} + 105\cdot2^{k+13+3i} + r_{8+i}(= 12285\cdot2^{19+4i})  & \text{if  } 0\leqslant i \leqslant l-4, \; k > 8+i \\
   3\cdot2^{2k+2l-1} + 137\cdot2^{k+3l+4} + r_{8+l-3}(=12029\cdot2^{23+4\cdot(l-4)})  & \text{if  } i = l-3,   \; k > 8+i  \\
 3\cdot2^{2k+2l+1} + 213\cdot2^{k+3l+7} + r_{8+l-2}(= 11373\cdot2^{27+4\cdot(l-4)})  & \text{if  } i = l-2,   \; k > 8+i \\
 11\cdot2^{2k+2l+3} + 333\cdot2^{k+3l+9} + r_{8+l-1} (= 10159\cdot2^{31+4\cdot(l-4)}) & \text{if  } i = l-1,  \; k > 8+i \\
 \end{cases}
\end{equation*}

   \begin{equation*}
\Gamma_{8+l+i}^{\left[2\atop{ 2+3\atop 2+3+l} \right]\times k} 
 = \begin{cases}
 21\cdot2^{2k+2l+5} + 2331\cdot2^{k+3l+11} + r_{8+l} (= 15435\cdot2^{34+4\cdot(l-4)}) & \text{if  } i = 0, \; k > 8+l+i  \\
    21\cdot2^{2k+2l+8} + 2163\cdot2^{k+3l+15} + r_{9+l}(= 2835\cdot2^{39+4\cdot(l-4)}) & \text{if  } i = 1, \; k > 8+l+i  \\
   53\cdot2^{2k+2l+11} + 1311\cdot2^{k+3l+19} + r_{10+l} (= -1417\cdot2^{44+4\cdot(l-4)} ) & \text{if  } i = 2,  \; k > 8+l+i \\
 \end{cases}
\end{equation*}    
   
    \begin{equation*}
\Gamma_{11+l+i}^{\left[2\atop{ 2+3\atop 2+3+l} \right]\times k} 
 = \begin{cases}
  105\cdot2^{2k+2l+14} -315\cdot2^{k+3l+23} + r_{11+l} (= 105\cdot2^{49+4\cdot(l-4)}) & \text{if  } i = 0,  \; k > 11+l+i \\
 2^{3k+l+9} - 7\cdot2^{2k+2l+18} +  7\cdot2^{k+3l+28} + r_{12+l} (= -2^{55+4\cdot(l-4)} ) &  \text{if   } i =1,\; k\geqslant 11+l+i    
\end{cases}
\end{equation*}  

Then we get:\\

   \begin{align*}
& R_{q} 
= 2^{(k+12+l)q -(3k+l+9)} \sum_{i = 0}^{\inf(k, 12+l)}
 \Gamma_{i}^{\left[2\atop {2+3\atop 2+3+l} \right]\times k} 2^{-iq}
   \end{align*}
 \end{example}  
 
 \subsection{The inverse problem in the case $s=2,\;m =3,\; l=4,\; k=6$}
\label{subsec 2}

 The results in this subsection are new.

    \begin{example} 
     We get:\\

  \begin{equation*}
\Gamma_{i}^{\left[2\atop{ 2+3\atop 2+3+4} \right]\times 6} 
 = \begin{cases}
 1  & \text{if  } i = 0, \\
 21  & \text{if  } i = 1, \\
 490 & \text{if  } i = 2, \\
 7200 & \text{if  } i = 3, \\
 108768 & \text{if  } i = 4, \\
 1679616 & \text{if  } i = 5, \\
2145687552 & \text{if  } i = 6, 
 \end{cases}
\end{equation*}  

 The number $ R_{q} $ of solutions \\
 $(Y_1,Z_1,U_{1},V_{1}, \ldots,Y_q,Z_q,U_{q},V_{q})  \in \big( \mathbb{F}_{2}[T ]  \big)^{4q} $ of the polynomial equations
   \[\left\{\begin{array}{c}
 Y_{1}Z_{1} +Y_{2}Z_{2}+ \ldots + Y_{q}Z_{q} = 0,  \\
   Y_{1}U_{1} + Y_{2}U_{2} + \ldots  + Y_{q}U_{q} = 0,\\
    Y_{1}V_{1} + Y_{2}V_{2} + \ldots  + Y_{q}V_{q} = 0,\\  
 \end{array}\right.\]
  satisfying the degree conditions \\
                   $$  degY_i \leq 5 , \quad degZ_i \leq 1 ,\quad degU_{i}\leq 4,\quad degV_{i}\leq 8 \quad for \quad 1\leq i \leq q. $$ \\    

is equal to\\

$ 2^{22q - 31} \sum_{i = 0}^{6}
 \Gamma_{i}^{\left[2\atop {2+3\atop 2+3+4} \right]\times 6} 2^{-iq} $
 \end{example}
 \newpage
 
  \subsection{The inverse problem in the case $s=3,\;m =0,\; l\geqslant 0$}
\label{subsec 3}

The results in this subsection are new.

  \begin{example}
\textbf{$s=3,\;m=0,\;l=0.$}
 
   The number  $ \Gamma_{i}^{\left[3\atop {3\atop 3} \right]\times k} $ of triple  persymmetric  $9\times k$ rank i   matrices over  $\mathbb{F}_{2}$   of the form \\
 $$   \left ( \begin{array} {cccccc}
\alpha _{1} & \alpha _{2}  &  \ldots  &  \ldots  & \alpha _{k-1}  &  \alpha _{k} \\
\alpha _{2 } & \alpha _{3} &  \ldots  &  \ldots &  \alpha _{k} &  \alpha _{k+1} \\
\alpha _{3 } & \alpha _{4} &   \ldots  &  \ldots  &  \alpha _{k+1} &  \alpha _{k+2} \\
\hline
\beta  _{1} & \beta  _{2}   &  \ldots  &  \ldots  &  \beta_{k-1} &  \beta _{k}  \\
\beta  _{2} & \beta  _{3} &  \ldots  &  \ldots  &    \beta_{k} &  \beta _{k+1}  \\
\beta  _{3} & \beta  _{4} &  \ldots  &  \ldots   &  \beta_{k+1} &  \beta _{k+2}  \\
\hline
\gamma  _{1} & \gamma   _{2} &  \ldots  &  \ldots  & \gamma  _{k-1} &  \gamma  _{k}  \\
\gamma  _{2} & \gamma  _{3}&  \ldots  &  \ldots    & \gamma  _{k} &  \gamma  _{k+1}  \\
\gamma  _{3} & \gamma   _{4} &  \ldots  &  \ldots    & \gamma  _{k+1} &  \gamma  _{k+2}  \\
 \end{array}  \right). $$ 
 is equal to

\[ \Gamma_{i}^{\left[3\atop{ 3\atop 3} \right]\times k}=
 \begin{cases}
1 &\text{if  }  i = 0 \\
21 &\text{if  }  i = 1 \\
378  &\text{if  }  i = 2 \\
7\cdot2^{k+ 2} +  5936  & \text{if   } i = 3 \\
147\cdot2^{k+ 2} + 84672  & \text{if   } i = 4 \\
147\cdot9\cdot2^{k+ 3} + 959616  & \text{if   } i = 5 \\
7\cdot2^{2k+4} + 2121\cdot2^{k+ 6} + 5863424   & \text{if   } i = 6\\
105\cdot2^{2k+4} + 2625\cdot2^{k+ 9} - 92897280  & \text{if   } i = 7 \\
105\cdot2^{2k+8} - 315\cdot2^{k+ 14} + 220200960  & \text{if   } i = 8 \\
2^{3k+6} - 7\cdot2^{2k+12} + 7\cdot2^{k+ 19} - 134217728  & \text{if  } i = 9,\;k\geq 9
\end{cases}
\]

\end{example}

 \begin{example}
\textbf{$s=3,\;m=0,\;l=1.$}

  The number  $ \Gamma_{i}^{\left[3\atop {3\atop 3+1} \right]\times k} $ of triple  persymmetric  $10\times k$ rank i   matrices over  $\mathbb{F}_{2}$   of the form \\
 $$   \left ( \begin{array} {cccccc}
\alpha _{1} & \alpha _{2}  &  \ldots  &  \ldots  & \alpha _{k-1}  &  \alpha _{k} \\
\alpha _{2 } & \alpha _{3} &  \ldots  &  \ldots &  \alpha _{k} &  \alpha _{k+1} \\
\alpha _{3 } & \alpha _{4} &   \ldots  &  \ldots  &  \alpha _{k+1} &  \alpha _{k+2} \\
\hline
\beta  _{1} & \beta  _{2}   &  \ldots  &  \ldots  &  \beta_{k-1} &  \beta _{k}  \\
\beta  _{2} & \beta  _{3} &  \ldots  &  \ldots  &    \beta_{k} &  \beta _{k+1}  \\
\beta  _{3} & \beta  _{4} &  \ldots  &  \ldots   &  \beta_{k+1} &  \beta _{k+2}  \\
\hline
\gamma  _{1} & \gamma   _{2} &  \ldots  &  \ldots  & \gamma  _{k-1} &  \gamma  _{k}  \\
\gamma  _{2} & \gamma  _{3}&  \ldots  &  \ldots    & \gamma  _{k} &  \gamma  _{k+1}  \\
\gamma  _{3} & \gamma   _{4} &  \ldots  &  \ldots    & \gamma  _{k+1} &  \gamma  _{k+2}  \\
\gamma  _{4} & \gamma  _{5} &  \ldots  &  \ldots  & \gamma  _{k+2} &  \gamma  _{k+3}  
 \end{array}  \right). $$ 
 is equal to

 \begin{equation*}
 \Gamma_{i}^{\left[3\atop {3\atop 3+1} \right]\times k}  
= \begin{cases}
 1  & \text{if  } i = 0, \\
 21  & \text{if  } i = 1, \\
 378  & \text{if  } i = 2, \\
3\cdot2^{k+2} +6192  & \text{if  } i = 3, \\
71\cdot2^{k+2} +124864  & \text{if  } i = 4, \\
651\cdot2^{k+3} +1246848  & \text{if  } i = 5, \\
2^{2k+4} + 645\cdot2^{k+7} +15464448  & \text{if  } i = 6, \\
27\cdot2^{2k+4} + 531\cdot2^{k+11} +93782016 & \text{if  } i = 7, \\
105\cdot2^{2k+6} + 2625\cdot2^{k+12} -1486356480 & \text{if  } i = 8, \\
105\cdot2^{2k+10} -315\cdot2^{k+17} +3523215360 & \text{if  } i = 9, \\
 2^{3k+7} - 7\cdot2^{2k+14} +  7\cdot2^{k+22} -2147483648 &  \text{if   } \; i =10      
\end{cases}
\end{equation*}  

\end{example}

 \begin{example}
\textbf{$s=3,\;m=0,\;l=2.$}

  The number  $ \Gamma_{i}^{\left[3\atop {3\atop 3+2} \right]\times k} $ of triple  persymmetric  $11\times k$ rank i  matrices over  $\mathbb{F}_{2}$  of the form \\
 $$   \left ( \begin{array} {cccccc}
\alpha _{1} & \alpha _{2}  &  \ldots  &  \ldots  & \alpha _{k-1}  &  \alpha _{k} \\
\alpha _{2 } & \alpha _{3} &  \ldots  &  \ldots &  \alpha _{k} &  \alpha _{k+1} \\
\alpha _{3 } & \alpha _{4} &   \ldots  &  \ldots  &  \alpha _{k+1} &  \alpha _{k+2} \\
\hline
\beta  _{1} & \beta  _{2}   &  \ldots  &  \ldots  &  \beta_{k-1} &  \beta _{k}  \\
\beta  _{2} & \beta  _{3} &  \ldots  &  \ldots  &    \beta_{k} &  \beta _{k+1}  \\
\beta  _{3} & \beta  _{4} &  \ldots  &  \ldots   &  \beta_{k+1} &  \beta _{k+2}  \\
\hline
\gamma  _{1} & \gamma   _{2} &  \ldots  &  \ldots  & \gamma  _{k-1} &  \gamma  _{k}  \\
\gamma  _{2} & \gamma  _{3}&  \ldots  &  \ldots    & \gamma  _{k} &  \gamma  _{k+1}  \\
\gamma  _{3} & \gamma   _{4} &  \ldots  &  \ldots    & \gamma  _{k+1} &  \gamma  _{k+2}  \\
\gamma  _{4} & \gamma  _{5} &  \ldots  &  \ldots  & \gamma  _{k+2} &  \gamma  _{k+3} \\
 \gamma  _{5} & \gamma  _{6} &  \ldots  &  \ldots  & \gamma  _{k+3} &  \gamma  _{k+4} \\

 \end{array}  \right). $$ 
 is equal to
 \begin{equation*}
   \Gamma_{i}^{\left[3\atop {3\atop 3+2} \right]\times k} 
= \begin{cases}
 1  & \text{if  } i = 0, \\
 21  & \text{if  } i = 1, \\
 378  & \text{if  } i = 2, \\
3\cdot2^{k+2} +6192  & \text{if  } i = 3, \\
39\cdot2^{k+2} +99264  & \text{if  } i = 4, \\
347\cdot2^{k+3} + 1503872  & \text{if  } i = 5, \\
2^{2k+4} + 618\cdot2^{k+6} +21426176  & \text{if  } i = 6, \\
3\cdot2^{2k+4} + 1293\cdot2^{k+9} +246448128  & \text{if  } i = 7, \\
27\cdot2^{2k+6} + 531\cdot2^{k+14} +1500512256 & \text{if  } i = 8, \\
105\cdot2^{2k+8} + 2625\cdot2^{k+15} - 92897280\cdot2^{8} & \text{if  } i = 9, \\
105\cdot2^{2k+12} -315\cdot2^{k+20} +53760\cdot2^{20} & \text{if  } i = 10, \\
 2^{3k+8} - 7\cdot2^{2k+16} +  7\cdot2^{k+25} - 2^{35} &  \text{if   }  i =11      
\end{cases}
\end{equation*}

\end{example}

 \begin{example}
\textbf{$s=3,\;m=0,\;l=3.$}

  The number  $ \Gamma_{i}^{\left[3\atop {3\atop 3+3} \right]\times k} $ of triple  persymmetric  $12\times k$ rank i  matrices over  $\mathbb{F}_{2}$  of the form \\
 $$   \left ( \begin{array} {cccccc}
\alpha _{1} & \alpha _{2}  &  \ldots  &  \ldots  & \alpha _{k-1}  &  \alpha _{k} \\
\alpha _{2 } & \alpha _{3} &  \ldots  &  \ldots &  \alpha _{k} &  \alpha _{k+1} \\
\alpha _{3 } & \alpha _{4} &   \ldots  &  \ldots  &  \alpha _{k+1} &  \alpha _{k+2} \\
\hline
\beta  _{1} & \beta  _{2}   &  \ldots  &  \ldots  &  \beta_{k-1} &  \beta _{k}  \\
\beta  _{2} & \beta  _{3} &  \ldots  &  \ldots  &    \beta_{k} &  \beta _{k+1}  \\
\beta  _{3} & \beta  _{4} &  \ldots  &  \ldots   &  \beta_{k+1} &  \beta _{k+2}  \\
\hline
\gamma  _{1} & \gamma   _{2} &  \ldots  &  \ldots  & \gamma  _{k-1} &  \gamma  _{k}  \\
\gamma  _{2} & \gamma  _{3}&  \ldots  &  \ldots    & \gamma  _{k} &  \gamma  _{k+1}  \\
\gamma  _{3} & \gamma   _{4} &  \ldots  &  \ldots    & \gamma  _{k+1} &  \gamma  _{k+2}  \\
\gamma  _{4} & \gamma  _{5} &  \ldots  &  \ldots  & \gamma  _{k+2} &  \gamma  _{k+3} \\
 \gamma  _{5} & \gamma  _{6} &  \ldots  &  \ldots  & \gamma  _{k+3} &  \gamma  _{k+4} \\
 \gamma  _{6} & \gamma  _{7} &  \ldots  &  \ldots  & \gamma  _{k+4} &  \gamma  _{k+5} 
 \end{array}  \right). $$ 
 is equal to
 \begin{equation*}
  \Gamma_{i}^{\left[3\atop {3\atop 3+3} \right]\times k} 
= \begin{cases}
 1  & \text{if  } i = 0, \\
 21  & \text{if  } i = 1, \\
 378  & \text{if  } i = 2, \\
3\cdot2^{k+2} +6192  & \text{if  } i = 3, \\
39\cdot2^{k+2} +99264  & \text{if  } i = 4, \\
219\cdot2^{k+3} + 1569408  & \text{if  } i = 5, \\
2^{2k+4} + 314\cdot2^{k+6} +24113152  & \text{if  } i = 6, \\
3\cdot2^{2k+4} + 621\cdot2^{k+9} +342786048 & \text{if  } i = 7, \\
3\cdot2^{2k+6} + 1293\cdot2^{k+12} +3943170048 & \text{if  } i = 8, \\
27\cdot2^{2k+8} + 531\cdot2^{k+17} + 93782016\cdot2^{8} & \text{if  } i = 9, \\
105\cdot2^{2k+10} + 2625\cdot2^{k+18} - 92897280\cdot2^{12} & \text{if  } i = 10, \\
105\cdot2^{2k+14} -315\cdot2^{k+23} +53760\cdot2^{24} & \text{if  } i = 11, \\
 2^{3k+9} - 7\cdot2^{2k+18} +  7\cdot2^{k+28} - 2^{39} &  \text{if   }  i =12      
\end{cases}
\end{equation*}

\end{example}

 \begin{example}
\textbf{$s=3,\;m=0,\;l=4.$}

  The number  $ \Gamma_{i}^{\left[3\atop {3\atop 3+4} \right]\times k} $ of triple  persymmetric  $13\times k$ rank i  matrices over  $\mathbb{F}_{2}$  of the form \\
 $$   \left ( \begin{array} {cccccc}
\alpha _{1} & \alpha _{2}  &  \ldots  &  \ldots  & \alpha _{k-1}  &  \alpha _{k} \\
\alpha _{2 } & \alpha _{3} &  \ldots  &  \ldots &  \alpha _{k} &  \alpha _{k+1} \\
\alpha _{3 } & \alpha _{4} &   \ldots  &  \ldots  &  \alpha _{k+1} &  \alpha _{k+2} \\
\hline
\beta  _{1} & \beta  _{2}   &  \ldots  &  \ldots  &  \beta_{k-1} &  \beta _{k}  \\
\beta  _{2} & \beta  _{3} &  \ldots  &  \ldots  &    \beta_{k} &  \beta _{k+1}  \\
\beta  _{3} & \beta  _{4} &  \ldots  &  \ldots   &  \beta_{k+1} &  \beta _{k+2}  \\
\hline
\gamma  _{1} & \gamma   _{2} &  \ldots  &  \ldots  & \gamma  _{k-1} &  \gamma  _{k}  \\
\gamma  _{2} & \gamma  _{3}&  \ldots  &  \ldots    & \gamma  _{k} &  \gamma  _{k+1}  \\
\gamma  _{3} & \gamma   _{4} &  \ldots  &  \ldots    & \gamma  _{k+1} &  \gamma  _{k+2}  \\
\gamma  _{4} & \gamma  _{5} &  \ldots  &  \ldots  & \gamma  _{k+2} &  \gamma  _{k+3} \\
 \gamma  _{5} & \gamma  _{6} &  \ldots  &  \ldots  & \gamma  _{k+3} &  \gamma  _{k+4} \\
 \gamma  _{6} & \gamma  _{7} &  \ldots  &  \ldots  & \gamma  _{k+4} &  \gamma  _{k+5} \\
  \gamma  _{7} & \gamma  _{8} &  \ldots  &  \ldots  & \gamma  _{k+5} &  \gamma  _{k+6} 
 \end{array}  \right). $$ 
 is equal to
 \begin{equation*}
  \Gamma_{i}^{\left[3\atop {3\atop 3+4} \right]\times k} 
=  \begin{cases}
 1  & \text{if  } i = 0, \\
 21  & \text{if  } i = 1, \\
 378  & \text{if  } i = 2, \\
3\cdot2^{k+2} +6192  & \text{if  } i = 3, \\
39\cdot2^{k+2} +99264  & \text{if  } i = 4, \\
219\cdot2^{k+3} + 1569408  & \text{if  } i = 5, \\
2^{2k+4} + 186\cdot2^{k+6} +25161728  & \text{if  } i = 6, \\
3\cdot2^{2k+4} + 317\cdot2^{k+9} +38577664 & \text{if  } i = 7, \\
3\cdot2^{2k+6} + 621\cdot2^{k+12} +5356032\cdot2^{10} & \text{if  } i = 8, \\
3\cdot2^{2k+8} + 1293\cdot2^{k+15} + 240672\cdot 2^{18} & \text{if  } i = 9, \\
27\cdot2^{2k+10} + 531\cdot2^{k+20} + 93782016\cdot 2^{12} & \text{if  } i = 10, \\
105\cdot2^{2k+12} + 2625\cdot2^{k+21} - 92897280\cdot 2^{16} & \text{if  } i = 11, \\
105\cdot2^{2k+16} -315\cdot2^{k+26} +53760\cdot 2^{28} & \text{if  } i = 12, \\
2^{3k+10} - 7\cdot2^{2k+20} +  7\cdot2^{k+31} - 2^{43} &  \text{if   }  i =13      
\end{cases}
\end{equation*}

\end{example}

 \begin{example}
\textbf{$s=3,\;m=0,\;l=5.$}

  The number  $ \Gamma_{i}^{\left[3\atop {3\atop 3+5} \right]\times k} $ of triple  persymmetric  $14\times k$ rank i  matrices over  $\mathbb{F}_{2}$  of the form \\
 $$   \left ( \begin{array} {cccccc}
\alpha _{1} & \alpha _{2}  &  \ldots  &  \ldots  & \alpha _{k-1}  &  \alpha _{k} \\
\alpha _{2 } & \alpha _{3} &  \ldots  &  \ldots &  \alpha _{k} &  \alpha _{k+1} \\
\alpha _{3 } & \alpha _{4} &   \ldots  &  \ldots  &  \alpha _{k+1} &  \alpha _{k+2} \\
\hline
\beta  _{1} & \beta  _{2}   &  \ldots  &  \ldots  &  \beta_{k-1} &  \beta _{k}  \\
\beta  _{2} & \beta  _{3} &  \ldots  &  \ldots  &    \beta_{k} &  \beta _{k+1}  \\
\beta  _{3} & \beta  _{4} &  \ldots  &  \ldots   &  \beta_{k+1} &  \beta _{k+2}  \\
\hline
\gamma  _{1} & \gamma   _{2} &  \ldots  &  \ldots  & \gamma  _{k-1} &  \gamma  _{k}  \\
\gamma  _{2} & \gamma  _{3}&  \ldots  &  \ldots    & \gamma  _{k} &  \gamma  _{k+1}  \\
\gamma  _{3} & \gamma   _{4} &  \ldots  &  \ldots    & \gamma  _{k+1} &  \gamma  _{k+2}  \\
\gamma  _{4} & \gamma  _{5} &  \ldots  &  \ldots  & \gamma  _{k+2} &  \gamma  _{k+3} \\
 \gamma  _{5} & \gamma  _{6} &  \ldots  &  \ldots  & \gamma  _{k+3} &  \gamma  _{k+4} \\
 \gamma  _{6} & \gamma  _{7} &  \ldots  &  \ldots  & \gamma  _{k+4} &  \gamma  _{k+5} \\
  \gamma  _{7} & \gamma  _{8} &  \ldots  &  \ldots  & \gamma  _{k+5} &  \gamma  _{k+6} \\
    \gamma  _{8} & \gamma  _{9} &  \ldots  &  \ldots  & \gamma  _{k+6} &  \gamma  _{k+7} 
   \end{array}  \right). $$ 
 is equal to
 \begin{equation*}
\Gamma_{i}^{\left[3\atop{ 3\atop 3+5} \right]\times k} 
 = \begin{cases}
 1  & \text{if  } i = 0, \\
 21  & \text{if  } i = 1, \\
 378  & \text{if  } i = 2, \\
3\cdot2^{k+2} +6192  & \text{if  } i = 3, \\
39\cdot2^{k+2} +99264  & \text{if  } i = 4, \\
219\cdot2^{k+3} + 1569408  & \text{if  } i = 5, \\
2^{2k+4} + 186\cdot2^{k+6} +25161728  & \text{if  } i = 6, \\
3\cdot2^{2k+4} + 189\cdot2^{k+9} +402554880 & \text{if  } i = 7, \\
3\cdot2^{2k+6} + 317\cdot2^{k+12} +385777664\cdot2^{4} & \text{if  } i = 8, \\
3\cdot2^{2k+8} + 621\cdot2^{k+15} +342786048\cdot2^{8} & \text{if  } i = 9, \\
3\cdot2^{2k+10} + 1293\cdot2^{k+18} + 246448128\cdot 2^{12} & \text{if  } i = 10, \\
27\cdot2^{2k+12} + 531\cdot2^{k+23} + 93782016\cdot 2^{16} & \text{if  } i = 11, \\
105\cdot2^{2k+14} + 2625\cdot2^{k+24} - 92897280\cdot 2^{20} & \text{if  } i = 12, \\
105\cdot2^{2k+18} -315\cdot2^{k+29} +53760\cdot 2^{32} & \text{if  } i = 13, \\
2^{3k+11} - 7\cdot2^{2k+22} +  7\cdot2^{k+34} - 2^{47} &  \text{if   }  i =14      
\end{cases}
\end{equation*}  
\end{example}

\begin{example}
\textbf{$s=3,\;m=0,\; l\geqslant 5$}

  The number  $ \Gamma_{i}^{\left[3\atop {3\atop 3+l} \right]\times k} $ of triple  persymmetric  $(9+l)\times k$ rank i  matrices over  $\mathbb{F}_{2}$  of the form \\
 $$   \left ( \begin{array} {cccccc}
\alpha _{1} & \alpha _{2}  &  \ldots  &  \ldots  & \alpha _{k-1}  &  \alpha _{k} \\
\alpha _{2 } & \alpha _{3} &  \ldots  &  \ldots &  \alpha _{k} &  \alpha _{k+1} \\
\alpha _{3 } & \alpha _{4} &   \ldots  &  \ldots  &  \alpha _{k+1} &  \alpha _{k+2} \\
\hline
\beta  _{1} & \beta  _{2}   &  \ldots  &  \ldots  &  \beta_{k-1} &  \beta _{k}  \\
\beta  _{2} & \beta  _{3} &  \ldots  &  \ldots  &    \beta_{k} &  \beta _{k+1}  \\
\beta  _{3} & \beta  _{4} &  \ldots  &  \ldots   &  \beta_{k+1} &  \beta _{k+2}  \\
\hline
\gamma  _{1} & \gamma   _{2} &  \ldots  &  \ldots  & \gamma  _{k-1} &  \gamma  _{k}  \\
\gamma  _{2} & \gamma  _{3}&  \ldots  &  \ldots    & \gamma  _{k} &  \gamma  _{k+1}  \\
\gamma  _{3} & \gamma   _{4} &  \ldots  &  \ldots    & \gamma  _{k+1} &  \gamma  _{k+2}  \\
\gamma  _{4} & \gamma  _{5} &  \ldots  &  \ldots  & \gamma  _{k+2} &  \gamma  _{k+3} \\
 \gamma  _{5} & \gamma  _{6} &  \ldots  &  \ldots  & \gamma  _{k+3} &  \gamma  _{k+4} \\
 \gamma  _{6} & \gamma  _{7} &  \ldots  &  \ldots  & \gamma  _{k+4} &  \gamma  _{k+5} \\
  \gamma  _{7} & \gamma  _{8} &  \ldots  &  \ldots  & \gamma  _{k+5} &  \gamma  _{k+6} \\
    \gamma  _{8} & \gamma  _{9} &  \ldots  &  \ldots  & \gamma  _{k+6} &  \gamma  _{k+7} \\
  \vdots  & \vdots  &   \vdots  & \vdots  &   \vdots  & \vdots                         \\
    \gamma  _{3+l} & \gamma  _{4+l} &  \ldots  &  \ldots  & \gamma  _{k+1+l} &  \gamma  _{k+2+l} \\
  \end{array}  \right). $$ 
 is equal to
 \begin{equation*}
\Gamma_{i}^{\left[3\atop{ 3\atop 3+l} \right]\times k} 
 = \begin{cases}
 1  & \text{if  } i = 0, \\
 21  & \text{if  } i = 1, \\
 378  & \text{if  } i = 2, \\
3\cdot2^{k+2} +6192  & \text{if  } i = 3, \\
39\cdot2^{k+2} +99264  & \text{if  } i = 4, \\
219\cdot2^{k+3} + 1569408  & \text{if  } i = 5, \\
2^{2k+4} + 186\cdot2^{k+6} +25161728  & \text{if  } i = 6, \\
3\cdot2^{2k+2i-10} + 189\cdot2^{k+3i-12} +402554880\cdot 2^{4(i-7)} & \text{if  } 7\leqslant i\leqslant l+2 \\
3\cdot2^{2k+2l-4} + 317\cdot2^{k+3l-3} +38577664\cdot2^{4l-16} & \text{if  } i = l+3, \\
3\cdot2^{2k+2l-2} + 621\cdot2^{k+3l} +342786048\cdot2^{4l-12} & \text{if  } i=l+4, \\
3\cdot2^{2k+2l} + 1293\cdot2^{k+3l+3} +246448128\cdot2^{4l-8} & \text{if  } i=l+5, \\
27\cdot2^{2k+2l+2} + 531\cdot2^{k+3l+8} +93782016\cdot2^{4l-4} & \text{if  } i=l+6, \\
105\cdot2^{2k+2l+4} + 2625\cdot2^{k+3l+9} -92897280\cdot2^{4l} & \text{if  } i=l+7, \\
105\cdot2^{2k+2l+8} -315\cdot2^{k+3l+14} +53760\cdot2^{4l+12} & \text{if  } i=l+8  \\
2^{3k+l+6} -7\cdot2^{2k+2l+12} +7\cdot2^{k+3l+19}-2^{4l+27} & \text{if  } i = l+9
\end{cases}
\end{equation*} 
\begin{proof}
We compute $ \Gamma_{i}^{\left[3\atop{ 3\atop 3+l} \right]\times k} $ for $0\leqslant i\leqslant l+9$   by using the following two reduction formulas:\\

\begin{align*} 
\Gamma_{7+j}^{\left[3\atop{ 3\atop 3+l} \right]\times k} 
 & =2^{4j}\cdot\Gamma_{7}^{\left[3\atop{ 3\atop 3+l-j} \right]\times (k-j)}  & \text{if  } 0\leq j\leq l,  \; k \geqslant 7+j   \\
 \Gamma_{7+l+j}^{\left[3\atop{ 3\atop 3+l} \right]\times k} 
 & =2^{4l+12j}\cdot\Gamma_{2(3-j)+1}^{\left[3-j\atop{ 3-j\atop 3-j} \right]\times (k-l-3j)}  & \text{if  } 0\leq j\leq 2,  \; k \geqslant 7+l+j 
\end{align*}

combined with the formula:\\
\begin{equation*}
\Gamma_{7}^{\left[3\atop{ 3\atop 3+l} \right]\times k} = \Gamma_{7}^{\left[3\atop{ 3\atop 3+5} \right]\times k} \quad   \text{if  }\; l\geqslant 5,  \; k >  7
\end{equation*}

More precisely : (with $ k > 7+j $)\\

 \begin{equation*}
\Gamma_{7+j}^{\left[3\atop{ 3\atop 3+l} \right]\times k} 
= \begin{cases}
2^{4j}\cdot\Gamma_{7}^{\left[3\atop{ 3\atop 3+5} \right]\times (k-j)} =
3\cdot2^{2k+2j+4} + 189\cdot2^{k+3j+9} +402554880\cdot2^{4j} & \text{if  } 0\leqslant j \leqslant l-5, \\
2^{4l-16}\cdot\Gamma_{7}^{\left[3\atop{ 3\atop 3+4} \right]\times (k-l+4)} =
3\cdot2^{2k+2l-4} + 317\cdot2^{k+3l-3} +38577664\cdot2^{4l-16} & \text{if  } j = l-4, \\
2^{4l-12}\cdot\Gamma_{7}^{\left[3\atop{ 3\atop 3+3} \right]\times (k-l+3)} =
3\cdot2^{2k+2l-2} + 621\cdot2^{k+3l} +342786048\cdot2^{4l-12} & \text{if  } j = l-3, \\
2^{4l-8}\cdot\Gamma_{7}^{\left[3\atop{ 3\atop 3+2} \right]\times (k-l+2)} =
3\cdot2^{2k+2l} + 1293\cdot2^{k+3l+3} +246448128\cdot2^{4l-8} & \text{if  } j = l-2, \\
2^{4l-4}\cdot\Gamma_{7}^{\left[3\atop{ 3\atop 3+1} \right]\times (k-l+1)} =
27\cdot2^{2k+2l+2} + 531\cdot2^{k+3l+8} +93782016\cdot2^{4l-4} & \text{if  } j = l-1, \\
2^{4l}\cdot\Gamma_{7}^{\left[3\atop{ 3\atop 3} \right]\times (k-l)} =
105\cdot2^{2k+2l+4} + 2625\cdot2^{k+3l+9} -92897280\cdot2^{4l} & \text{if  } j = l, \\
2^{4l+12}\cdot\Gamma_{5}^{\left[2\atop{ 2\atop 2} \right]\times (k-l-3)} =
105\cdot2^{2k+2l+8} -315\cdot2^{k+3l+14} +53760\cdot2^{4l+12} & \text{if  } j = l+1, \\
2^{4l+24}\cdot\Gamma_{3}^{\left[1\atop{ 1\atop 1} \right]\times (k-l-6)} =
2^{3k+l+6} -7\cdot2^{2k+2l+12} +7\cdot2^{k+3l+19}-2^{4l+27} & \text{if  } j = l+2
\end{cases}
\end{equation*} 
\end{proof}
 The number $ R_{q} $ of solutions \\
 $(Y_1,Z_1,U_{1},V_{1}, \ldots,Y_q,Z_q,U_{q},V_{q})  \in \big( \mathbb{F}_{2}[T ]  \big)^{4q} $ of the polynomial equations
   \[\left\{\begin{array}{c}
 Y_{1}Z_{1} +Y_{2}Z_{2}+ \ldots + Y_{q}Z_{q} = 0,  \\
   Y_{1}U_{1} + Y_{2}U_{2} + \ldots  + Y_{q}U_{q} = 0,\\
    Y_{1}V_{1} + Y_{2}V_{2} + \ldots  + Y_{q}V_{q} = 0,\\  
 \end{array}\right.\]
  satisfying the degree conditions \\
                   $$  degY_j \leq k-1 , \quad degZ_j \leq 2 ,\quad degU_{j}\leq 2,\quad degV_{j}\leq 2+l \quad for \quad 1\leq j \leq q, $$ \\    
   is equal to              
   \begin{align*}
2^{(k+9+l)q -(3k+l+6)} \sum_{i = 0}^{\inf(k, 9+l)}
 \Gamma_{i}^{\left[3\atop {3\atop 3+l} \right]\times k} 2^{-iq}
   \end{align*}
In particular:  \\

   \begin{align*}
 R_{1}= 2^{-2k+3} \sum_{i = 0}^{\inf(k, 9+l)}
 \Gamma_{i}^{\left[3\atop {3\atop 3+l} \right]\times k} 2^{-i} = 2^{9+l} +2^{k} -1
   \end{align*}

\end{example}

 \selectlanguage{english}

  \end{document}